\documentclass{amsart}

\usepackage[mathcal]{euscript}
\usepackage{amsfonts,amssymb,amsbsy,amsmath}
\usepackage{graphicx,color}
\usepackage{amscd}

\def\Dg:{\endgraf{\bf Dg:}\enspace\ignorespaces}
\def\n:{\endgraf{\bf n:}\enspace\ignorespaces}
\def\new:{\endgraf{\bf new:}\enspace\ignorespaces}

\def\latin#1{\emph{#1}}
\def\ie{\latin{i.e.}}
\def\cf.{\latin{cf\.}}

\def\frak{\mathfrak}
\def\Cal{\mathcal}
\def\Bbb{\mathbb}
\def\bold{\mathbf}

\def\onto{\twoheadrightarrow}
\def\into{\hookrightarrow}
\let\<\langle
\let\>\rangle
\let\MG\Gamma
\def\tMG{\tilde\MG}
\def\CH{\Cal H}

\def\ba{\bold a}
\def\bb{\bold b}
\def\BG#1{\Bbb B_{#1}}
\def\FG#1{\Bbb F_{#1}}
\def\SG#1{\Bbb S_{#1}}
\def\CG#1{{\Z}_{#1}}
\def\DG#1{\Bbb D_{#1}}
\def\diagram{\Cal D}
\def\necklace{\Cal N}

\def\cpic#1{$\vcenter{\hbox{%
\input{\GLEfdirectory#1.inc}}}$}

\let\Gb\beta
\let\Gg\gamma
\let\Gd\delta
\let\Gf\varphi
\let\Gs\sigma

\def\1{^{-1}}
\def\C{\Bbb C}
\def\R{\Bbb R}
\def\X{\Bbb X}
\def\Y{\Bbb Y}
\def\Z{\Bbb Z}
\def\Cp#1{\Bbb{P}^{#1}}
\def\Rp#1{\Cp{#1}_\R}
\def\bm{\frak m}
\def\tbm{\tilde{\bm}}
\def\BM{\bar{\bm}}

\def\bminf{\bm_{\infty}}
\def\trans{^{\mathrm{t}}}
\let\sminus\smallsetminus

\makeatletter
\def\ra@{\mathbin{\scriptstyle\uparrow}}
\def\la@{\mathbin{\scriptstyle\downarrow}}
\def\ra{\protect\ra@}
\def\la{\protect\la@}

\def\rom#1{\leavevmode
  \edef\prevskip@{\ifdim\lastskip=\z@ \else\hskip\the\lastskip\relax\fi}%
  \unskip
  \edef\prevpenalty@{\ifnum\lastpenalty=\z@ \else
    \penalty\the\lastpenalty\relax\fi}%
  \unpenalty \/\prevpenalty@ \prevskip@ {\rm #1}}

\makeatletter
\def\:{\nobreak\hskip.1111em\mathpunct{}\nonscript\mkern-\thinmuskip{:}\hskip
 .3333emplus.0555em\relax}
\let\DOTaccent\.
\def\.{.\spacefactor\@m}
\makeatother

\theoremstyle{plain}
\newtheorem{t.}{Theorem}[section]
\newtheorem{l.}[t.]{Lemma}
\newtheorem{p.}[t.]{Proposition}
\newtheorem{c.}[t.]{Corollary}

\makeatletter
\expandafter\let\expandafter\c@equation\csname c@t.\endcsname
\expandafter\let\expandafter\theequation\csname thet.\endcsname
\makeatother

\theoremstyle{definition}

\newtheorem{d.}[t.]{Definition}

\newtheorem{r.}[t.]{Remark}
\newtheorem*{w.}{Warning}
\newtheorem*{conv.}{Convention}

\def\theorem{\begin{t.}}
\def\endtheorem{\end{t.}}
\def\lemma{\begin{l.}}
\def\endlemma{\end{l.}}
\def\proposition{\begin{p.}}
\def\endproposition{\end{p.}}
\def\corollary{\begin{c.}}
\def\endcorollary{\end{c.}}
\def\definition{\begin{d.}}
\def\enddefinition{\end{d.}}
\def\remark{\begin{r.}}
\def\endremark{\end{r.}}
\def\Figure{\begin{figure}}
\def\endFigure{\end{figure}}
\def\Table{\begin{table}}
\def\endTable{\end{table}}

\def\done{\qed}
\def\pni{\qed}
\def\donesymbol{$\vartriangleleft$}
\def\done{{\def\qedsymbol{\donesymbol}\qed}\rm}
\def\pnisymbol{$\vartriangleright$}
\def\pni{{\def\qedsymbol{\pnisymbol}\qed}\rm}

\makeatletter

\def\roster@*{\begin{itemize}\parskip\z@\def\endroster{\end{itemize}}}
\newcommand\roster@@[1][\@ne]{\begin{enumerate}\parskip\z@
 \setcounter{enumi}#1
 \global\advance\c@enumi\m@ne
 \def\endroster{\end{enumerate}}}
\def\roster{\@ifnextchar*\roster@\roster@@}

\def\equation@*{\begin{equation*}\def\]{\end{equation*}}}
\def\equation@@{\begin{equation}\def\]{\end{equation}}}
\def\[{\@ifnextchar*\equation@\equation@@}

\makeatother

\def\epsfdirectory{gle/}
\newif\ifincludefiles
\includefilestrue

\def\cpic#1{$\vcenter{\hbox{\ifincludefiles

 \input{\epsfdirectory#1.inc}\else

 \includegraphics{\epsfdirectory#1}%
\fi}}$}

\makeatletter
\def\tabstrut{\vrule height10\p@ depth3\p@ width\z@}
\def\toprule{\noalign{\hrule}}
\def\midrule{\noalign{\hrule}}
\def\bottomrule{\noalign{\hrule}}
\makeatother

\def\hide#1{\rlap{\smash{$#1$}}}

\makeatletter
\def\Circ{\raise\p@\hbox{$\m@th\scriptstyle\bigcirc$}}
\def\Square{\raise-.7\p@\hbox{$\m@th\square$}}
\def\ndrel{{\joinrel\relbar\joinrel}}
\def\ndstone#1{\setboxz@h{$\joinrel\relbar\joinrel\joinrel\relbar\joinrel$}%
 \hbox to\wdz@{\hss$#1$\hss}\kern-\wdz@\boxz@}
\def\ND#1{\ND@#1\ND@}
\def\ND@#1{\let\next\ND@\ndrel
 \ifx<#1\ndstone<\else
 \ifx>#1\ndstone>\else
 \ifx.#1\Circ\else
 \ifx*#1\Square\else
 \ifx\ND@#1\let\next\relax\else
 #1\fi\fi\fi\fi\fi\next}
\def\nCirc{\NO{\scriptscriptstyle\!\bigcirc}}
\def\nSquare{\NO\square}

\def\Mup{{\uparrow}}
\def\Mdown{{\downarrow}}
\def\Mside{{*}}
\def\Mskip@{\mskip1mu}
\def\Mcurve#1{\Mcurve@#1\Mcurve@}
\def\Mcurve@#1{\Mskip@\let\next\Mcurve@
 \ifx u#1\Mup\else
 \ifx d#1\Mdown\else
 \ifx.#1\Mside\else
 \ifx\Mcurve@#1\let\next\Mskip@\else
 #1\fi\fi\fi\fi\next}

\def\inserthyphen{\ifcat\next a-\fi\ignorespaces}
\let\BLACK\bullet
\let\WHITE\circ
\def\SMALLER#1#2{\vcenter{\hbox{$\m@th#1\mathord#2$}}}
\def\CROSS{\SMALLER\scriptstyle\times}
\let\STAR*
\def\TRIANG{\SMALLER\scriptstyle\vartriangle}
\def\pblack-{$\BLACK$\futurelet\next\inserthyphen}
\def\pwhite-{$\WHITE$\futurelet\next\inserthyphen}
\def\pcross-{$\CROSS$\futurelet\next\inserthyphen}
\def\pstar-{$\STAR$\futurelet\next\inserthyphen}
\def\ptriang-{$\TRIANG$\futurelet\next\inserthyphen}
\def\black{\protect\pblack}
\def\white{\protect\pwhite}
\def\cross{\protect\pcross}

\def\triang{\protect\ptriang}
\def\NO#1{\mathord\#_{#1}}
\def\nblack{\NO\BLACK}
\def\nwhite{\NO\WHITE}

\def\ncirc{\NO\diamond}

\def\graph{\frak G}
\makeatother

\def\ls|#1|{\mathopen|#1\mathclose|}

\makeatletter

\def\dessin{\frak D}
\def\skeleton{\frak S}
\def\bm{\frak{m}}
\def\BM{\bar\bm}
\let\disk\Omega

\def\I{{\operator@font I}}
\def\II{{\operator@font II}}
\def\III{{\operator@font III}}
\def\IV{{\operator@font IV}}

\def\qopname@@#1{\mathop{\operator@font#1}\nolimits}
\def\qopname@{\protect\qopname@@}
\let\opname\qopname@

\def\Aut{\qopname@{Aut}}
\def\Fix{\qopname@{Fix}}
\def\Dssn{\qopname@{Dssn}}
\def\conv{\qopname@{conv}}
\def\val{\qopname@{val}}
\def\Sk{\qopname@{Sk}}
\def\Im{\qopname@{Im}}
\def\Stab{\qopname@{Stab}}
\def\stab{\qopname@{stab}}
\def\id{\qopname@{id}}
\def\Map{\qopname@{Map}}
\def\trace{\qopname@{trace}}

\def\group@@#1{\qopname@@{\hbox{\sl#1}}}
\def\group@{\protect\group@@}
\def\SL{\group@{SL}}
\def\GL{\group@{GL}}
\def\PSL{\group@{PSL}}
\def\PGL{\group@{PGL}}
\def\O{\group@{O}}
\def\SO{\group@{SO}}
\def\NS{\group@{NS}}
\def\MW{\group@{MW}}
\def\BMG{\group@{MG}}
\def\BND{\group@{BND}}

\makeatother

\begin{document}

\title[Products of pairs of Dehn twists]
{Products of pairs of Dehn twists\\ and maximal real Lefschetz fibrations}

\author{Alex Degtyarev}
\address{Bilkent University\\
Department of Mathematics\\
06800 Ankara, Turkey}
\email{degt@fen.bilkent.edu.tr}

\author{Nerm{\accent95\i}n Salepc{\accent95\i}}
\address{
Institut Camille Jordan\\
Universit\'{e} Lyon I\\
43, Boulevard du 11 Novembre 1918\\
69622 Villeurbanne Cedex, France}
\email{salepci@math.univ-lyon1.fr}


\thanks{%
The second author
was partially supported by
the European
Community's Seventh Framework Programme ([FP7/2007-2013] [FP7/2007-2011])
under grant agreement no\. [258204],
as well as
by the French \emph{Agence nationale de la recherche} grant
ANR-08-BLAN-0291-02}

\begin{abstract}
We address the problem of existence and uniqueness of a  factorization of a
given element of the modular group into a product of two Dehn twists. As a
geometric application, we conclude that any maximal real elliptic Lefschetz
fibration is algebraic.
\end{abstract}

\maketitle

\section{Introduction}

\subsection{Motivation}
An object repeatedly occurring  in algebraic geometry is a fibration with
singular fibers. If the base is a topological disk $D^2$ and the number of
singular fibers is finite, the topology (and, in some extremal cases, the
analytic structure as well) can adequately be described  by the so-called
\emph{monodromy factorization} of the monodromy at infinity (the boundary
of~$D^2$).

More precisely, consider a
proper smooth map
$p\:X\to B\cong D^2$ and  let
$\Delta:=\{b_{1},b_{2},\dots,b_{r}\}$
be the set of the critical values of
$p$, which are all assumed in the interior of $B$.
The restriction of $p$ to $B^\sharp:=B\smallsetminus \Delta$ is a locally
trivial fibration and one can consider its \emph{monodromy}
$\bm\:\pi_{1}(B^\sharp, b)\to \Aut F_{b}$, where $F_{b}$ is the fiber over a
fixed base point $b\in B^\sharp$ and $G:=\Aut F_{b}$
 is an appropriately defined group of classes of automorphisms of $F_{b}$.
 (The precise nature of
the automorphisms used and their equivalence depend on a particular problem.)
The monodromy at infinity
$\bminf:=\bm[\partial B]\in G$
is usually assumed fixed in advance.

\begin{w.}
Throughout
the paper, all group actions are \emph{right}.
(It is under this convention that monodromy is a homomorphism.) This convention
applies to matrix groups as well: our matrices act on row vectors by the right
multiplication.
Given a right action $X\times G\to X$, we
denote by
$x\ra g$
the image of $x\in X$ under $g\in G$.
\end{w.}

Consider a system of lassoes, one lasso $\gamma_{i}$ about each critical
value $b_{i}$, $i=1,\ldots, r$,  disjoint except at the common base point~$b$
and  such that $\gamma_{1}\cdot \ldots \cdot \gamma_{r}\sim \partial B$. (Such
a system is called a \emph{geometric basis} for $\pi_{1}(B^\sharp,b)$.)
Evaluating the monodromy~$\bm$ at each~$\gamma_i$, we
obtain a sequence $\bm_{i}:=\bm(\gamma_{i})$.

\begin{d.}
Given a group $G$, a $G$-valued \emph{monodromy factorization} of length $r$
is a finite ordered
sequence $\BM:=(\bm_{1},\ldots,\bm_{r})$ of elements of~$G$.
The product
$\bminf:=\bm_{1}\cdot \ldots \cdot\bm_{r}$ is called the \emph{monodromy at
infinity }of $\BM$, and $\BM$ itself is often referred to as a monodromy
factorization of $\bminf$. The subgroup of $G$ generated by
$\bm_{1},\ldots,\bm_{r}$ is called the \emph{monodromy group} of $\BM$.
\end{d.}

The ambiguity in the choice of a geometric basis leads to a certain
equivalence relation. According to Artin~\cite{Artin}, if $b\in \partial B$,
any two geometric bases are related by an element of the braid group $\BG{r}$.
Hence, the corresponding monodromy factorizations are related by a sequence
of \emph{Hurwitz moves} \[
\sigma_{i}\:(\ldots,\bm_{i},\bm_{i+1},\ldots) \mapsto
(\ldots,\bm_{i}\bm_{i+1}\bm_{i}^{-1},\bm_{i},\ldots),\quad i=1,\ldots, r-1.
\label{eq.Hurwitz}
\]
If the base point is not on the boundary or if the identification between
$F_{b}$
and
the `standard'
fiber is not fixed,  one should also consider
the \emph{global conjugation}
\[*
g^{-1}\BM g=(g^{-1}\bm_{1}g,\ldots,g^{-1}\bm_{r}g)
\]
by an element $g\in G$.

\begin{d.}
Two monodromy factorizations are said to be \emph{strongly} (\emph{weakly})
Hurwitz equivalent if they can be related by a finite sequence of Hurwitz
moves (respectively, a sequence of Hurwitz moves and  global conjugation). For
brevity, we routinely simplify this term to just
strong/weak
equivalence.
\end{d.}

It is immediate that both the monodromy at infinity and the monodromy group
are invariant under strong Hurwitz equivalence, whereas their conjugacy
classes are invariant under weak
Hurwitz equivalence.

The most well known examples where this machinery applies are
\roster*
\item
ramified coverings, with $G=\SG{n}$ the symmetric group;
\item
algebraic or, more generally, pseudo holomorphic curves in $\C^2$, with
$G=\BG{n}$ the braid group;
\item
(real) elliptic surfaces or, more generally,
(real) genus one
Lefschetz fibrations, with $G=\tMG:=\SL(2,\Z)$
 the mapping class group of a
torus.
\endroster
(Literature on the subject is abundant, and we direct the reader
to~\cite{degt:monodromy} for further references.)
Typically, the topological type of
a singular fiber $F_{i}:=p^{-1}(b_{i})$
is determined by the conjugacy class of the corresponding element $\bm_{i}$,
and it is common to restrict the topological types by assuming that all
$\bm_{i}$ should belong to a certain
preselected set of conjugacy classes.
Thus,   in the three examples above,
`simplest'
singular fibers correspond to,
respectively, transpositions in $\SG{n}$, Artin generators in $\BG{n}$, and
Dehn twists in~$\tMG$, see subsection~\ref{ss.prem}.

A monodromy factorization satisfying this additional restriction is often
called \emph{simple}, and
a
wide open problem with a great deal of possible
geometric implications is the classification, up to
strong/weak Hurwitz
equivalence, of the simple monodromy factorizations of
a given element $\bminf\in G$
and
of a given length.

\subsection{Principal results}

Geometrically, our principal subject is elliptic Lefschetz fibrations, and
the algebraic counterpart is the classification of the factorizations of a
given element $\bminf \in \tMG$ into products of Dehn twists. At this point,
it is worth mentioning that there are
cyclic
central extensions $ \tMG \onto \MG$
and  $\BG{3}\onto \MG$,  where $\MG:=\PSL(2,\Z)$
is
the modular group, and each
Dehn twist in $\MG$ lifts to a unique
Dehn twist in $\tMG$ or, respectively, to a
unique Artin generator in $\BG{3}$; hence, the problems of the
classification of
\emph{simple}  monodromy factorizations in all three groups are equivalent.
For this reason, we will mainly work in $\MG$. To simplify the further
exposition, we introduce the following terminology: an
\emph{$r$-factorization}
(of an element $g\in\Gamma$) is  a monodromy factorization
$\BM=(\bm_{1},\ldots,\bm_{r})$ with each $\bm_{i}$ a Dehn twist and such that
$\bminf=g$. To shorten the notation, we will often speak about an
$r$-factorization
$g=\bm_{1}\cdot \ldots \cdot \bm_{r}$.

Even with the group as simple as $\BG{3}$ (the first non-abelian braid
group), surprisingly little is known. On the one hand, according to
Moishezon--Livn\'{e}~\cite{Moishezon:LNM}, a
$6k$-factorization
of a power
$(\sigma_{1}\sigma_{2})^{3k}$
of the Garside element is unique up to strong Hurwitz equivalence.
This result was recently generalized by Orevkov~\cite{Orevkov:talk} to any
element positive in the standard Artin basis $\sigma_{1}, \sigma_{2}$. On the
other hand,
a series of exponentially large (in $r$) sets of non-equivalent
$r$-factorizations of the same element $g_r:=L^{5r-6}\in\MG$
(depending on~$r$) was recently constructed
in~\cite{degt:monodromy}; furthermore,
these factorizations are
indistinguishable by most conventional invariants.
(For some other examples,
related to the next braid group~$\BG{4}$,
see~\cite{Kulikov.Auroux.Shevchishin}.)

Thus, it appears that, in its full generality, the problem of the
classification of the $r$-factorizations of a given element is rather
difficult and quite far from its complete understanding. In this paper, we
confine ourselves to $2$-factorizations only, addressing both their existence
and uniqueness. Even in this simplest case, the results obtained seem rather
unexpected.

Algebraically, our principal results are the three theorems below. For
the statements, we briefly recall that the elements of the modular group are
commonly divided into elliptic, parabolic, and hyperbolic, the former being those
of finite order, and the two latter being those that, up to conjugation, can
be represented by a word in positive powers of a particular pair $L, R$ of
generators of $\Gamma$, see subsection~\ref{ss.prem} for further details.
(Whenever speaking about words in a given alphabet, we mean \emph{positive}
words only; if negative powers are allowed, they are listed in the alphabet
explicitly.)
We
use $A^{\mathrm{t}}$ for the transpose of a matrix  $A$.
One has $L\trans=R$; hence, the transpose $A\trans$ of a word~$A$ in $\{L,R\}$ is again a word
in $\{L,R\}$: it is obtained from~$A$ by interchanging
$L\leftrightarrow R$ and reversing the order of the letters.

\theorem\label{th.>=1}
An element $g\in\MG$ admits a $2$-factorization if and only if either
\roster
\item
$g\sim\X=RL\1$ \rom($g$ is elliptic\rom), or
\item
$g\sim R^2$ or $g\sim L^4$ \rom($g$ is parabolic\rom), or
\item
$g\sim L^2AL^2A\trans$ for some word~$A\ne\varnothing$ in $\{L,R\}$
\rom($g$ is hyperbolic\rom).
\endroster
\endtheorem

\theorem\label{th.<=2}
The number of weak equivalence classes
of  $2$-factorizations of  $g\in\MG$ is at most one if $g$ is elliptic or parabolic, and at
most two if $g$ is hyperbolic.
\endtheorem

\theorem\label{th.weak=strong}
The single weak equivalence class of  $2$-factorizations of an element
$g\sim L^4$ splits into two strong equivalence
classes\rom:
\[*
L^4=R\cdot(R\1L^2)R(R\1L^2)\1
=LRL\1\cdot(LR\1L^2)R(LR\1L^2)\1.
\]
In all
other cases, each weak equivalence class of  $2$-factorizations
constitutes a single
strong equivalence class.
\endtheorem

Theorems~\ref{th.>=1}, \ref{th.<=2}, \ref{th.weak=strong} are proved in
subsections~\ref{proof.>=1}, \ref{proof.<=2}, \ref{proof.weak=strong},
respectively. The proofs are based on a relation between subgroups of the
modular group and a certain class of Grothendieck's  \emph{dessins d'enfants}.
A refinement of Theorem~ \ref{th.<=2}, namely a detailed description of the
elements admitting more than one $2$-factorization, is found in
subsection~\ref{s.=2}, see Theorem~\ref{th.=2}.

Another interesting  phenomenon related to the modular group is the fact that
some of its elements are \emph{real}, \ie,
they can be represented as a product of
two involutive elements of $\PGL(2,\Z)\smallsetminus \MG$. (For a geometric
interpretation
and further details, see \cite{Salepci:real}
and subsection~\ref{ss.real}.)  The
relation between this property and the existence/uniqueness of a
$2$-factorization,
as well as the existence of real $2$-factorizations, are
discussed in Theorem~\ref{th.2.and.real}.

Geometrically, $2$-factorizations are related to real relatively minimal
Jacobian elliptic Lefschetz fibrations over the sphere $S^2$ with two pairs
of complex conjugate singular fibers; an important class of such
fibrations are
some
maximal ones.
Intuitively, an elliptic Lefschetz
fibration is a topological counterpart of an algebraic elliptic surface
(see
section~\ref{S.Lf} for the precise definitions), and
one of the major questions
is the realizability of a given
\emph{real} elliptic
Lefschetz fibration by an algebraic one. (In the complex case, the answer to
this question is trivially in the affirmative due to the classification found
in \cite{Moishezon:LNM}, see Theorem~\ref{th.Moishezon.Livne};
in the real case, examples of non-algebraic fibrations are known,
see~\cite{Salepci:necklace,Salepci:thesis}.)
A
real Lefschetz fibration is \emph{maximal} if its real part has the maximal
Betti number with respect to the Thom--Smith
inequality~\eqref{eq.Thom-Smith}.
A maximal
real Lefschetz fibration may have $0$, $1$ or $2$ pairs of
complex conjugate singular
fibers, see \ref{s.counts}.
In the former case, the fibration is called \emph{totally real}, and
such a fibration is necessarily algebraic due to the following theorem.

\theorem[see~\cite{Salepci:necklace,Salepci:thesis}]\label{th.Nermin.alg}
Any totally real maximal Jacobian Lefschetz fibration
is algebraic.
\pni
\endtheorem

Amongst the most important
geometric applications
of the algebraic results of
the
paper is an extension of
Theorem~\ref{th.Nermin.alg} to all maximal Jacobian fibrations.

\theorem\label{th.algebraic}
Any maximal
Jacobian Lefschetz fibration is algebraic.
\endtheorem

This theorem is proved in subsection~\ref{s.algebraic}.

As another geometric application, we
settle
a question left unanswered in
\cite{DIK:elliptic}. Namely, we show that the equivariant deformation class
of a nonsingular real trigonal $M$-curve in a Hirzebruch surface (see
section~\ref{S.tc} for the definitions) is determined by the topology of its real
structure, see Theorems~\ref{th.def=dif}
and~\ref{th.def=dif.z=0}. Moreover, at
most two such curves may share homeomorphic real parts.

One may speculate that it is the relation to \emph{maximal} geometric
objects, which are commonly known to be topologically `rigid',
that makes $2$-factorizations relatively `tame'. At present, we do not have
any clue on what the general statements concerning the existence and
uniqueness of $r$-factorizations may look like. One of the major reasons is
the fact that,
even though an analogue of Proposition~\ref{44} holds for any number of Dehn
twists,
Lemma~\ref{interesting.lemma} does not have
a literate
extension to free groups on more than two generators,
\cf.~\cite{Bardakov}.

To our knowledge,
even the finiteness of the number of equivalence classes of factorizations of
a given element is still an open question. According to R.~Matveyev and K.~Rafi (private
communication), certain finiteness statements do hold in hyperbolic groups;
alas, neither~$\MG$ nor~$\BG3$ is hyperbolic. On the other hand, found in
B.~Moishezon~\cite{Moishezon:infty} is an example of an infinite sequence of
non-equivalent factorizations (although non-simple) of the element $\Delta^2$
in the braid group~$\BG{54}$.

\subsection{Contents of the paper}
Sections~\ref{S.mg}, \ref{S.Lf}, and \ref{S.tc} are of an auxiliary nature:
we recall the basic notions and necessary known results concerning,
respectively, the modular group, (real) elliptic Lefschetz fibrations, and
(real) trigonal curves. The heart of the paper is Section~\ref{S.heart}: the
principal algebraic results and their refinements are proved here.
Section~\ref{S.appl} deals with the geometric applications: we establish the
semi-simplicity of real trigonal $M$-curves and, as an upshot, prove
Theorem~\ref{th.algebraic}.

We use the conventional symbol \qedsymbol\ to mark the ends of the proofs.
Some statements are marked with \donesymbol\ or \pnisymbol: the former means
that the proof has already been explained (for example, most corollaries),
and the latter indicates that the proof is not found in the paper and the
reader is directed to the literature, usually cited at the beginning of the
statement.

\subsection{Acknowledgment}
This paper
was
essentially completed during the second author's stay as a
Leibniz fellow and the first author's visit as a \emph{Forschungsgast} to
the \emph{Mathematisches Forschungsinstitut Oberwolfach}; we are grateful  to
this institution and its friendly staff for their hospitality and for the
excellent working conditions. We would
like to thank
Viatcheslav
Kharlamov for
his encouragement and  interest in the subject,
and Alexander Klyachko, who brought to our attention the Frobenius type formulas
counting solutions
to equations in finite groups.
We
are also grateful
to Anton Klyachko and
to the anonymous referee of this text, who drew our
attention to
Bardakov's paper~\cite{Bardakov} and
Kulkarni's paper~\cite{Kulkarni}, respectively.

\section{The modular group}\label{S.mg}

\subsection{Presentations of $\Gamma$}\label{ss.prem}

Consider
$\CH=\Z\ba\oplus\Z\bb$, a rank two free abelian group with the skew-symmetric
bilinear form  $\bigwedge^2\CH\to\Z$  given by  $\ba\cdot\bb=1$.  We regard
$\tMG :=\SL(2,\Z)$ as a group acting on $\CH$. Moreover, $\tMG$ is the
group of symplectic auto-symmetries of $\CH$; it is generated by the
matrices
\[*
\X= \bmatrix1 &-1\\1&\phantom{-}0\endbmatrix,\quad
\Y =\bmatrix \phantom{-}0 &1\\-1& 0\endbmatrix
\]
such that $\X^3=-\id$, $\Y^2=-\id$.

The \emph{modular group} $\Gamma := \PSL(2,\Z)$ is the quotient
$\SL(2,\Z)/\pm\id$.  When it does note lead to a confusion, we use the same
notation for a matrix $A$ in $\tMG$ and its projection to $\MG$.
It is known that $\MG\cong \Z_{3}*\Z_{2}$; we will work with the
following two presentations of this group:
\[*
\Gamma =\left<\X, \Y : \X^3=\Y^2=\id \right>=
\left<L,R : RL\1R=L\1RL\1,\ (RL\1)^3=\id\right>,
\]
where
\[*
L= \bmatrix1 & 1\\0&1\endbmatrix= \X\Y ,\quad
R = \bmatrix1 &0\\1&1\endbmatrix=\X^2\Y,
\]
so that $\X=RL\1$ and $\Y=LR\1L=R\1LR\1$ in $\MG$.
For future references note that the powers of these matrices are given by
\[*L^{n} = \bmatrix1 & n\\0&1\endbmatrix,\quad
R^{n} = \bmatrix1&0\\n&1\endbmatrix,\quad
n\in\Z.
\]

Since $\Gamma$ is a free product of cyclic groups, we have the following statement.

\lemma\label{lem.centralizer}
Two elements
$f,g\in\MG$ commute if and only if they generate a cyclic subgroup, or,
equivalently, if they are both powers of a common element
$h\in\MG$.
\done
\endlemma

\subsection{The conjugacy classes}\label{ss.conjclas}
\begin{figure}[h]
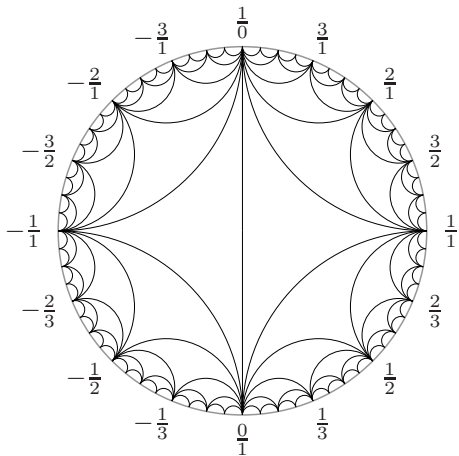

\cpic{Farey}
\caption{Poincar\'{e} disk endowed with the Farey tessellation}\label{farey}
\end{figure}

A
simple way to understand  the conjugacy classes is \latin{via}  the action of
$\Gamma$ on the Poincar\'{e} disk. The group $\Gamma$  is known to be the
symmetry group of  the Poincar\'{e} disk endowed with the
so-called
Farey tessellation,
shown in Figure~\ref{farey}. 
The non trivial elements of $\Gamma$ form three basic  families,
\emph{elliptic}, \emph{parabolic}, and \emph{hyperbolic}. These families are
distinguished by the nature of their fixed points on the Poincar\'{e} disk, or
equivalently, by the absolute value of their traces. Namely, an elliptic
matrix has $\mathopen|\trace\mathclose|<2$, so that it has a single fixed
point in the interior of the Poincar\'{e} disk and acts as a rotation with
respect to this fixed point. Elliptic matrices are the only torsion elements
of
$\Gamma$. A parabolic matrix has $\mathopen|\trace\mathclose|=2$; it has a
single rational fixed point (on the boundary of the Poincar\'{e} disk) and acts
as a rotation fixing this boundary point. A hyperbolic matrix, defined
\latin{via} $\mathopen|\trace\mathclose|>2$, has two irrational fixed points
on the boundary and acts as a translation fixing the geodesic connecting
these fixed points.

There are three conjugacy classes of elliptic matrices. Representatives of these classes can be taken as:
\[*
\Y=\bmatrix \phantom-0 & 1\\-1&0\endbmatrix,\quad
\X=\bmatrix 1& -1\\ 1&\phantom-0\endbmatrix,\quad
\X^{-1}=\bmatrix \phantom-0& 1\\ -1&1\endbmatrix.
\]

An element in $\Gamma$ is called a (positive) \emph{Dehn twist} if it is
conjugate to $R$ (the geometric meaning of this definition is explained in
subsection~\ref{ss.MCG}).
Any parabolic element is conjugate to
a certain
$n^{\mathrm{th}}$ power of a Dehn twist. Thus, a representative of a class
can be taken as $R^n$.

\begin{w.}
For the experts,
we emphasize that,
in accordance with our right group action convention, it is~$R$,
not~$L$, that represents a positive Dehn twist.
\end{w.}

The conjugacy classes of hyperbolic elements of~$\MG$ are determined by
sequences $[a_{1}, a_{2},\ldots,a_{2n}], a_{i}\in\Z^+$, defined up to
\emph{even} permutations and called \emph{cutting period cycles}.
Indeed, the fixed points of a hyperbolic matrix are irrational points that
are the zeroes of a quadratic equation, and they have a continued fraction expansion with the periodic tale
\[*
\ldots\, a_1+\dfrac1{a_2+\dfrac1{\,\,\ddots{\atop\dfrac1{a_{2n}}.}}}
\]
Note that $[a_{1}, a_{2},\ldots,a_{2n}]$ is not necessarily the minimal period: all matrices sharing the same pair of eigenvectors are powers of a minimal one, and the precise multiple of the minimal period corresponding to a given matrix $A$ can be recovered from its trace.

A representative of the conjugacy class corresponding to a cutting period cycle $[a_{1}, a_{2},\ldots,a_{2n}]$ can be chosen as
\[*
R^{a_{1}}\cdot L^{a_{2}}\cdot\ldots \cdot L^{a_{2n}}=
\bmatrix1 & 0\\a_{1}&1\endbmatrix\cdot \bmatrix1 & a_{2}\\0&1\endbmatrix\cdot \ldots \cdot \bmatrix1 & a_{n}\\0&1\endbmatrix.
\]
In the sequel, we will  be interested not only in the cutting period cycle
but also in the underlying word, called the \emph{cutting word}, in
two letters
$\{L,R\}$.
Recall that the cutting word encodes
the two types (right/left) of  triangles of the Farey tessellation cut by the
invariant
 geodesic, \cf.
\cite{Series,Salepci:real}. In terms of the cutting word, hyperbolic
conjugacy classes
can be
characterized as
those represented
by a word in
$\{L,R\}$
\emph{with both $L$ and $R$ present}.
Since
the cutting word is only defined up to cyclic
permutation, it is convenient to represent it
in the unit circle,
placing the
letters constituting the word
at equal angles (\cf. Figure~\ref{fig.n2} on page~\pageref{fig.n2}).
The resulting circle marked with a number of copies of~$L$
and~$R$ is called the \emph{cyclic diagram} $\diagram_g$ of a hyperbolic element~$g$.
One can also speak about the cyclic diagram of a parabolic element, with the
letters either all~$R$
(for a positive power of a Dehn twist)
or all~$L$
(for a negative power).

\subsection{Real elements}\label{ss.real} An involutive element of the coset $\PGL(2,\Z) \smallsetminus \Gamma$ is called a \emph{real structure} on $\Gamma$.
An element of $\MG$ is called \emph{real} if, in $\PGL(2,\Z)$,  it has a
decomposition into a product of  two real structures. For any real structure
$\tau$, let us define an
involutive
anti-automorphism $\hat{\tau}\:\Gamma \to \Gamma$
given by $\hat{\tau}(g)=\tau g^{-1}\tau$.
Then, a real element can also be defined as
one fixed by
$\hat{\tau}$ for some real
structure $\tau$. The significance of real elements is in their  geometric
interpretation. For example, such an element appears as the $\Gamma$-valued
monodromy at infinity of  a real elliptic Lefschetz fibration over a disk.

The characterization of real elements in $\Gamma$, as well as in $\tMG$, is
known, see \cite{Salepci:real}: all elliptic and parabolic matrices are
real, and a hyperbolic matrix is real if and only if its cutting period cycle
is \emph{odd bipalindromic}, \ie, up to cyclic permutation, it is a union of
two palindromic pieces of odd length. This property can be interpreted in
terms of the cyclic diagram as  the existence of a symmetry axis such that
the diagram is invariant under the reflection with respect to this axis.

Up
to conjugation, there are exactly two real structures on $\MG$:
\[
{\tau}_{1}=\bmatrix 0 & 1 \\1&0\endbmatrix,\quad
{\tau}_{2}=\bmatrix1 & \phantom-0\\0&-1\endbmatrix.\label{eq.realstruc}
\]
In the rest of the paper, $\tau_{1}$ and $\tau_{2}$ refer to these particular
matrices.
The action of $\hat\tau_i$ on the generators is as follows: \[
\begin{gathered}
\hat{\tau}_{1}(L)=R\1,\quad
\hat{\tau}_{1}(R)=L\1,\quad
\hat{\tau}_{1}(\X)=\X,\quad
\hat{\tau}_{1}(\Y)=\Y,\\
\hat{\tau}_{2}(L)=L,\quad
\hat{\tau}_{2}(R)=R,\quad
\hat{\tau}_{2}(\X)=\Y\X\Y,\quad
\hat{\tau}_{2}(\Y)=\Y.
\end{gathered}
\]

We extend the anti-automorphism $\hat{\tau}\:\MG\to\MG$ to the set of
$\MG$-valued monodromy factorizations as follows:
\[ \hat{\tau}(\bm_{1},\ldots, \bm_{r})=(\hat{\tau}(\bm_{r}),\ldots,\hat{\tau}(\bm_{1})).\]
(Note the reverse order.) It is straightforward that
the factorizations
$\hat{\tau}(\BM')$ and $\hat{\tau}(\BM'')$ are strongly/weakly equivalent if and only if so are $\BM'$ and $\BM''$. Furthermore, one has  $\hat{\tau}(\BM)_{\infty}=\hat{\tau}(\bminf)$ and the monodromy group of $\hat{\tau}(\BM)$ is the image of that of  $\BM$ under $\hat{\tau}$.

With future applications in mind, we will also discuss real
$2$-factorizations. A $2$-factorization
$\BM$
is said to be
\emph{real} if there is a real structure $\tau$ such that either
$\hat\tau(\BM)=\BM\ra\Gs_1$, see~\eqref{eq.Hurwitz}, or $\hat\tau(\BM)=\BM$.
The monodromy at infinity of a real $2$-factorization is obviously real; the
converse is not true, see \cite{Ozturk.Salepci:Giroux, Salepci:real}
and subsection~\ref{ss.relationreal}.

\remark\label{rem.bm.real.geometric}
Geometrically, a real
$2$-factorization represents a real Jacobian Lefschetz fibration over the
unit disk $D^2\subset\C$ (with the standard real structure $z\mapsto\bar z$)
with two singular fibers,
see subsection~\ref{s.Lf.real};
in the former case
($\hat\tau(\BM)=\BM\ra\Gs_1$),
the two singular fibers are real; in the latter case ($\hat\tau(\BM)=\BM$), they are complex
conjugate.
A specific example of a non-real $2$-factorization
with real monodromy at infinity is studied in~\cite{Ozturk.Salepci:Giroux};
this example has interesting geometric implications.
\endremark

\remark\label{rem.bm.real}
Alternatively, a $2$-factorization~$\BM$ is real if and only if $\hat\tau(\BM)$ is strongly
Hurwitz equivalent to~$\BM$ for some real structure~$\tau$. Indeed, since
an even power $\Gs_1^{2k}$ acts \latin{via} the conjugation
by the $\tau$-real element
$\bminf^{-k}$, it can be `undone' by replacing~$\tau$ with
$\tau':=\tau\bminf^k$, which is also a real structure.
In particular, it follows that being real is a property of a whole strong
Hurwitz equivalence class.
\endremark

\subsection{The mapping class group}\label{ss.MCG}
The \emph{mapping class group} $\Map_{+}(S)$ of an oriented smooth surface
$S$ is defined as the group of isotopy classes of orientation preserving
diffeomorphisms of $S$. If $S$
is the $2$-torus $T^2$, one can fix an
isomorphism $H_{1}(T^2,\Z)\cong \CH=\Z\ba\oplus\Z\bb$, and the map
$f\mapsto f_{*}$ establishes an isomorphism $\Map_{+}(T^2)\to \tMG$.

The (positive) \emph{Dehn twist }along a simple closed curve $l\in S$ is a
diffeomorphism of $S$ obtained by cutting $S$ along $l$ and regluing with a
twist of $2\pi$. If $S\cong T^2$, the image of the Dehn twist in the mapping
class group $\tMG$ depends only on the homology class $u :=[l]\in \CH$
and is given by the symplectic reflection
$x\mapsto x+(u,x)u$, where $(u, x)$ denotes the algebraic sum of the points of intersection of $u$ an $x$;
we denote
this image by $t_{u}$ and call it a \emph{Dehn twist} in $\tMG$. All  Dehn twists
form a whole conjugacy class which contains $R$;  they project to the
positive Dehn twists in $\MG$ introduced in subsection~\ref{ss.conjclas}.

\subsection{Subgroups of $\Gamma$}\label{ss.subgroups}
In this section, we summarize
the relation between the subgroups of~$\MG$ and a special class of bipartite
ribbon graphs, which we call \emph{skeletons}. A similar approach, in terms
of special triangulations of surfaces, was developed in~\cite{Bogomolov}. Our
approach is identical to the \emph{bipartite cuboid graphs}
in~\cite{Kulkarni}, except that we are mainly interested in subgroups of
infinite index and therefore are forced to consider infinite graphs supported
by non-compact surfaces.
We only recall briefly the few definitions and facts needed in the sequel;
for details and further references,
see~\cite{degt:monodromy}.
Note that, due to our right group action convention,
some definitions given below differ
slightly from those in~\cite{degt:monodromy}.

Recall that a \emph{ribbon graph} is a graph
(locally finite {\sl CW\/}-complex of dimension one), possibly infinite,
equipped
with a cyclic order
(\latin{i.q\.} transitive $\Z$-action)
on the star of each vertex. Typically, a ribbon graph is
a graph embedded into an oriented surface $S$, and the cyclic order is induced
by the orientation of $S$. In fact, a ribbon graph $\graph$ defines a unique,
up to homeomorphism, minimal oriented surface $S_{0}$
(non-compact if $\graph$ is infinite)
into which it is
embedded. The connected components of the complement
$S_{0}\smallsetminus\graph$ are called the \emph{regions} of $\graph$.

A  \emph{bipartite graph} is a graph whose vertices are colored with two colors:
\black-, \white-,  so that each edge connects vertices of opposite colors.

\definition\label{def.skeleton}
A
\emph{skeleton} is a connected bipartite ribbon graph with all
\black-vertices of valency~$3$ or~$1$ and all \white-vertices of valency~$2$
or~$1$. A skeleton is \emph{regular} if all its \black-- and \white-vertices
have valency~$3$ and~$2$, respectively.
\enddefinition

Since
$\MG=\{\X, \Y: \X^3=\Y^2=\id\}$, the set of edges of any skeleton is a
transitive $\MG$-set, with the action of~$\X$ and~$\Y$ given by the
distinguished cyclic order on the stars of, respectively, \black-- and
\white-vertices. (Due to the valency restrictions in
Definition~\ref{def.skeleton}, this action of $\Z*\Z$ does factor
through~$\MG$.)
Conversely, any transitive $\MG$-set can be regarded as (the
set of edges of) a skeleton, the \black-- and \white-vertices being the
orbits of~$\X$ and~$\Y$, respectively.
In the sequel, we identify the two categories.

As
a consequence, to each subgroup $G\subset\MG$ one can associate the
skeleton $G\backslash\MG$ (the set of left $G$-cosets, regarded as a right
$\MG$-set). This skeleton is regular if and only if $G$ is torsion free, \ie,
contains no elliptic elements; in this case, $G$ is free.
The skeleton $G\backslash \Gamma$ is equipped with a distinguished edge
$e:=G\backslash G$, which we call the \emph{base point}.
Conversely, given a skeleton $\skeleton$
and a base point $e$, the stabilizer $G$ of $e$ is a subgroup of $\Gamma$,
and one has $\skeleton=G\backslash \Gamma$. In general, without a base point
chosen, the \emph{stabilizer} of $\skeleton$ is defined as a conjugacy class
of subgroups of $\Gamma$.

\begin{conv.}
In the figures, we
usually
omit
most
bivalent \white-vertices, assuming that such a vertex is to be inserted
at the center of each `edge' connecting a pair of \black-vertices.
When of interest, the base point is denoted by a grey diamond.
For infinite skeletons, only a compact part is drawn and each maximal Farey
branch, see subsection~\ref{SS.pseudo.tree} and Figure~\ref{looptree}, left,
below, is represented by a \triang-vertex.
\end{conv.}

A combinatorial path (called a \emph{chain} in \cite{degt:monodromy})
in a skeleton~$\skeleton$ can be
regarded as  a pair $\gamma:=(e',g)$, where $e' $ is an edge, called the
\emph{initial point} of $\gamma$, and $g\in \Gamma$. Then $e'':=e'\ra g$ is
the \emph{terminal point of $\gamma$},
and the \emph{evaluation map}
$\val\:\Gg\mapsto g$ sends a path $\gamma=(e',g)$ to
its underlying element $g\in \Gamma$.
For a regular skeleton $\skeleton$,
the map $\val$ establishes
an isomorphism $\pi_{1}(\skeleton,e)=G$.
(In the presence of monovalent vertices, one should replace $\pi_1$ with an
appropriate orbifold fundamental group.)
When the initial point is understood,
we identify a path
$\Gg$ and its image $\val\Gg\in\MG$.
The \emph{product} of two paths is defined as usual:
$(e',g')\cdot(e'',g'')=(e',g'g'')$ provided that
$e''=e'\ra g'$;
the
\emph{inverse} of $\Gg=(e',g)$ is $\Gg\1:=(e'\ra g,g\1)$.

In the case of skeletons,
a \emph{region} can be redefined as an orbit of
$L=\X\Y$.
In this definition, a region~$\Cal R$
is the set of edges in the boundary of the
geometric realization of~$\Cal R$ whose canonical orientation
$\BLACK{\joinrel\rightarrow\joinrel\joinrel\relbar\joinrel}\WHITE$
agrees with the boundary orientation;
the other edges in the boundary are of the form $e\ra\Y$, $e\in\Cal R$.
An
\emph{$n$-gonal} region is an orbit of length $n$; intuitively, $n$ is the
number of \black-vertices in the boundary.
The minimal supporting surface~$S_0$ of a skeleton~$\skeleton$ can be
obtained by patching the boundary of each region~$\Cal R$ with a disk
(if $\Cal R$ is finite) or half-plane (if $\Cal R$ is infinite).

Given
a subgroup $G\subset\MG$,
the $G$-conjugacy classes of the Dehn twist
contained
in $G$ are
in a canonical one-to-one correspondence with the monogonal regions of the
skeleton $G\backslash\MG$, see~\cite{degt:monodromy}:
under the canonical identification
$G=\pi_1^{\text{orb}}(G\backslash\MG,G\backslash G)$ described above, these
classes are realized by the boundaries of the monogons.

\subsection{Pseudo-trees}\label{SS.pseudo.tree}
\begin{figure}[h]
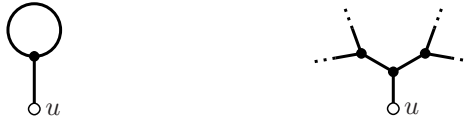

\cpic{01-loop}\hspace{3cm}
\cpic{01-tree}
\caption{A simple loop and a Farey branch}\label{looptree}
\end{figure}

\begin{figure}[h]
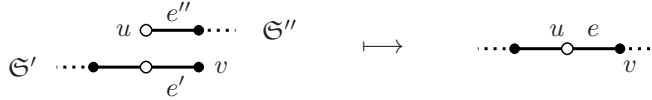

\cpic{splice}
\caption{The splice of two skeletons}\label{fig.splice}
\end{figure}

A special class of skeletons can be obtained from ribbon trees as follows.
Consider a ribbon tree with all \black-vertices of valency 3 (nodes) or 1
(leaves) and take its bipartite subdivision,
\ie, divide each edge into two
by inserting an extra \white-vertex in the middle.
We denote the resulting
graph by $\graph$. Let  us consider a \emph{vertex function} $\ell\:
\{\mbox{leaves}\}  \to \{0, \TRIANG, \BLACK, \WHITE\}$ such that, if two
leaves are incident to a common node, then $\ell$ does not assign \triang-\
to both. We perform the following modifications
at each leaf~$v$ of~$\graph$:
\roster*
\item
if $\ell(v)=\BLACK$, then no modification is done;
\item
if $\ell(v)=\WHITE$, then cut out  the leaf and the incident edge,
so that the resulting graph have a monovalent \white-vertex;
\item
if $\ell(v)=0$, then splice $\graph$ with a \emph{simple loop}, see
Figure~\ref{looptree}, left;
\item
if $\ell(v)=\TRIANG$, then splice $\graph$ with a \emph{Farey branch},
see Figure~\ref{looptree}, right.
\endroster
Formally,
a simple loop
is the skeleton
$\Gamma_{1}(2)\backslash \Gamma$, where
$\Gamma_1(2)=\left\{ \bmatrix1 & 0 \\ * &1\endbmatrix \bmod 2\right\}$,
and a Farey branch
$\Y \backslash \Gamma$
is
the only bipartite ribbon tree regular except
a single monovalent vertex, which is \white-.
Given two skeletons $\skeleton'$, $\skeleton''$, a monovalent
\black-vertex~$v$ of~$\skeleton'$, and a monovalent \white-vertex~$u$
of~$\skeleton''$, the \emph{splice} is defined as the skeleton obtained from
the disjoint union $\skeleton'\sqcup\skeleton''$ by identifying the
edges~$e'$, $e''$
incident to~$v$, $u$, respectively, to a common edge~$e$,
see Figure~\ref{fig.splice}.

A skeleton that can be obtained by the above procedure
is called a \emph{pseudo-tree}.
A pseudo-tree is regular if and only if
the images of $\ell$ are in $\{0, \TRIANG\}$.

Crucial for the sequel is the following statement, which
is an immediate consequence
from \cite[Proposition~4.4]{degt:monodromy}.

\proposition\label{44}
A proper subgroup $G\subset\Gamma$ is generated by two
distinct Dehn twists if and only if
its skeleton
$\skeleton:=G\backslash\MG$
is a regular pseudo-tree
with exactly two simple loops.
In this case, $G$ is freely generated by two Dehn twists. \pni
\endproposition

\begin{figure}[h]
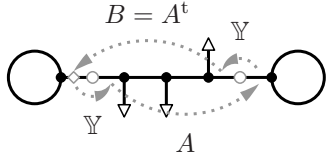

\cpic{01-pseudo}
\caption{An example of a pseudo-tree}\label{pseudo}
\end{figure}

Due to the requirement on the
\triang-values of a vertex function, a pseudo-tree $\skeleton$ as
in Proposition~\ref{44} looks as shown in Figure~\ref{pseudo}.
More
precisely, $\skeleton$ consists of two monogons connected by a horizontal
line segment and a number of Farey branches, upward and downward, attached
to this segment.
Thus, starting from one of the monogons, one can encode~$\skeleton$ and,
hence, the subgroup~$G$ itself by the sequence of the directions (up/down) of
the Farey branches.

\remark\label{rem.AB}
The \emph{monodromy at infinity} of
a
pseudo-tree $\skeleton$ is the conjugacy class $\bminf$ in $\Gamma$
realized by a large circle encompassing the compact part of $\skeleton$. Let
us choose the base point $e$ next to one of the monogons as shown in
Figure~\ref{pseudo}. Starting from $e$, we can realize $\bminf$ by the element
$L^2AL^2B$, where $A$ and $B$ are the paths shown in the figure.  Namely,
$A$ starts at  $e':=e\ra (\X\Y)^2=e\ra L^2=e\ra\Y$
and is a product of copies of
$R=\X^2\Y$ and $L=\X\Y$,  each downward \triang-vertex contributing an $R$
and each upward  \triang-vertex contributing an $L$. The other path $B$ can
be described similarly starting from a base point next to the other monogonal
region. However, it is obvious from the figure that the loop
$(e, \Y A \Y B)$
is contractible. Hence,  $B=\Y A^{-1} \Y$, and one can easily verify that
$A\trans=\Y A^{-1} \Y$ in $\Gamma$. Thus, we arrive at
\[\bminf\sim L^2 A L^2 A\trans,\label{m.infty}\]
where the word~$A$ in $\{L,R\}$ (possibly empty) is as described above.
As an upshot of this description we have the converse statement: a representation of the monodromy at infinity in the form~\eqref{m.infty} determines a pseudo-tree up to isomorphism.
\endremark

\section{The classification of $2$-factorizations}\label{S.heart}

\subsection{Proof of Theorem~\ref{th.>=1}}  \label{proof.>=1}
We precede the proof of this theorem with a few auxiliary statements.

\lemma\label{lem.tu.tv}
Two Dehn twists $t_u$, $t_v$, $u,v\in\CH$,
generate~$\MG$ if and only if $u$ and~$v$ span~$\CH$. If this is the case, the pair $(t_u,t_v)$ is conjugate to $(R,L\1)$.
\endlemma

\begin{proof}
If $u$ and~$v$ span~$\CH$, the signs can be chosen so that the
matrix~$M$ formed by $u$, $v$ as rows has determinant~$1$, \ie, belongs to~$\tMG$. The
conjugation by~$M$ takes $(R,L\1)$ to $(t_u,t_v)$; hence, $t_u$
and~$t_v$ generate~$\MG$.

For the converse statement, assume that the subgroup $\CH'\subset\CH$
spanned by~$u$ and~$v$ is proper. Since Dehn twists are
symplectic
reflections, see subsection~\ref{ss.MCG},
the subgroup~$\CH'$ is
obviously invariant under the subgroup $\MG'\subset\MG$ generated
by~$t_u$ and~$t_v$. Thus, there are primitive vectors in~$\CH$
that are not in the orbit $u\ra\MG'$. On the other hand, all
primitive vectors are known to form a single $\MG$-orbit; hence,
$\MG'\subset\MG$ is a proper subgroup.
\end{proof}

\lemma[\cf. Bardakov~\cite{Bardakov}]\label{interesting.lemma}
Let
$G:=\<\alpha,\beta\>$ be a free group, and let $\alpha',\beta'\in G$ be two
elements generating~$G$ and such that each~$\alpha'$, $\beta'$ is
conjugate to one of the original generators~$\alpha$, $\beta$. Then the
pair $(\alpha',\beta')$ is weakly Hurwitz equivalent to~$(\alpha,\beta)$.
\endlemma

\begin{proof}
After a global conjugation, possibly followed by~$\Gs_1$, one can
assume that $\alpha'=\alpha$. Then obviously $\beta'=T\1\beta T$ for some reduced
word~$T$ in $\{\alpha^{\pm1},\beta^{\pm1}\}$. One can assume that the first
letter of~$T$ is not~$\beta^{\pm1}$ and, after a global conjugation by
a power of~$\alpha$, one can also assume that the last letter of~$T$ is
not~$\alpha^{\pm1}$. Then, after expressing~$\alpha'$ and~$\beta'$ in terms
of~$\alpha$ and~$\beta$, any reduced word in $\{(\alpha')^{\pm1},(\beta')^{\pm1}\}$
results in a reduced word: no cancelation occurs. On the other
hand, there is a word that is equal to~$\beta$. Hence, one must have
$T=\id$ and $\beta'=\beta$.
\end{proof}

\begin{proof} [Proof of Theorem~\ref{th.>=1}]
Let $g\in \Gamma$
be
an element together with a $2$-factorization $\BM=(\bm_{1},\bm_{2})$.
Denote by $G$ the monodromy group of $\BM$.

If  $G$ is $\Gamma$, then by  Lemma~\ref{lem.tu.tv}
the pair $(\bm_{1},\bm_{2})$ is conjugate to $(R, L^{-1})$, and thus  $g$ is conjugate to
$\X=R\cdot L^{-1}$, which is an elliptic element.

If $\bm_{1}=\bm_{2}$, then $G$ is a cyclic subgroup of $\Gamma$;
hence, $g$ is conjugate to  $R^2=R\cdot R$, which is a parabolic element.

Otherwise, by Proposition~\ref{44},  $G$ is a proper subgroup such that
$G\backslash \Gamma$ is a regular pseudo-tree $\skeleton$ with two simple
loops. On $\skeleton$, choose a base point $e$ next to one of the monogons
and fix generators $\alpha, \beta$ of $G=\pi_{1}(\skeleton,e)$ as
shown in
Figure~\ref{markpoint}.
\begin{figure}[h]
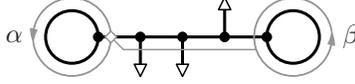

\cpic{01-pseudo2}
\caption{Generators of $G$ with respect to the base point}\label{markpoint}
\end{figure}
By Lemma~\ref{interesting.lemma}, the pair $(\alpha,
\beta)$ is weakly Hurwitz equivalent to $(\bm_{1},\bm_{2})$. Therefore, we
get $g\sim \bminf\sim L^2 AL^2 A\trans$, see \eqref{m.infty}.
If $A=\varnothing$, we get a parabolic element $g\sim L^4$; all other
elements obtained in this way are hyperbolic.

To finish the proof, note that the three cases mentioned above give the
complete list of subgroups generated by two Dehn twists,
and the conditions listed in the statement are necessary.
For the sufficiency, observe that a
factorization $g\sim L^2AL^2A\trans$ is not only
a necessary condition
but also a description of a particular $2$-factorization,
with the two Dehn twists as follows:\
\begin{align}
\label{eq.decomposition}
L^2AL^2A\trans&=
(\X L^{-1} \X^{-1} \Y)
(A)
(\X L^{-1} \X^{-1} \Y)
(\Y A^{-1} \Y)\\&=
\X L\1\X\1\cdot(\Y A\X)L\1(\Y A\X)\1.\qedhere
\end{align}
\end{proof}

Although the
converse statements are
contained in the above discussions, let us underline
the relation between the type of an element and the monodromy group of its
$2$-factorization.

\corollary[of the proof]\label{cor.type->group}
The monodromy group $G$ of any $2$-factorization of an element $g\in \Gamma$
is as follows\rom:
\roster*
\item
$g\sim \X$ \rom(elliptic\rom)  if and only if $G=\Gamma$\rom;
\item
$g\sim R^2$ \rom(parabolic\rom) if and only if
$G\subset\MG$ is a cyclic subgroup
generated by a single Dehn twist\rom;
\item
$g\sim L^4$ \rom(parabolic\rom) or $g$ is hyperbolic if and only if
$G\subset\MG$ is
a subgroup
as in Proposition~\ref{44}.
\done
\endroster
\endcorollary

\begin{figure}[h]
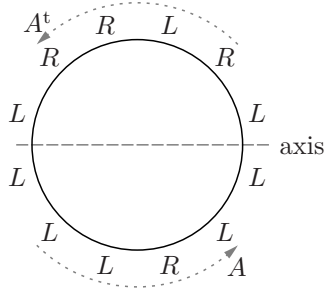

\cpic{n-2}
\caption{Cyclic diagram associated to $L^2AL^2A^t$ and its para-symmetry}\label{fig.n2}
\end{figure}

\remark\label{rem.para}
Geometrically,
a representation of an element $g$
in the form \eqref{m.infty} and the factorization \eqref{eq.decomposition}
can be described in terms of a \emph{para-symmetry} on the cyclic diagram of
$g$.  Let us call
the four special copies of~$L$ in
the word $L^2AL^2A\trans$
\emph{anchors}.
On the cyclic diagram,
trace an axis
passing
between the two anchors constituting each of the two pairs~$L^2$,
see Figure~\ref{fig.n2}. The reflection
with respect to this axis preserves the four anchors, while reversing the types
of all other letters. A reflection with this properties is called a
\emph{para-symmetry}.
We underline that the anchors are always of type~$L$.
\endremark

\corollary[of the proof and Remark~\ref{rem.AB}]\label{cor.A1A2}
The $2$-factorizations~\eqref{eq.decomposition} resulting from two
representations $L^2A_{1}L^2A_{1}\trans\sim L^2A_{2}L^2A_{2}\trans$ of the
same conjugacy class
are weakly
equivalent if and only if $A_{1}=A_{2}$ or $A_{1}=A_{2}\trans$. \done
\endcorollary

\subsection{Proof of Theorem~\ref{th.<=2}}\label{proof.<=2}
If $g$ is an elliptic element, we can assume that $g=\X=R\cdot L\1$.
Given another $2$-factorization $\X=t_u\cdot t_v$,
the two Dehn twists must generate~$\MG$, see Corollary~\ref{cor.type->group}.
Then,
due to Lemma~\ref{lem.tu.tv}, we have $t_u=h\1Rh$ and $t_v=h\1L\1h$ for some
$h\in\MG$. It follows that $h$ centralizes~$\X$ and hence $h$ is a power
of~$\X$, see Lemma~\ref{lem.centralizer}; thus, the second $2$-factorization
is strongly equivalent to the first one (as the
conjugation by the
monodromy at infinity
is the Hurwitz move $\Gs_1^{-2}$).

The only $2$-factorization of the parabolic element $g=R^2$ is $R^2$ itself,
as two distinct Dehn twists would produce either~$\X$, or~$L^4$, or a
hyperbolic element,
see
Corollary~\ref{cor.type->group}.
Finally, the parabolic element $g\sim L^4$ can be regarded as $V_0$,
see~\eqref{eq.Vm}, and this case
is considered below.
The two
orthogonal para-symmetries of the cyclic diagram of~$g$ result in two
conjugate (by~$L$) $2$-factorizations, which are not strongly
equivalent, as the corresponding marked skeletons (\cf. Figure~\ref{fig.L4}
on page~\pageref{fig.L4})
are not
isomorphic, see subsection~\ref{ss.subgroups}.

Now, assume that $g$ is a hyperbolic element and consider its cyclic
diagram $\diagram:=\diagram_g$. By assumption, it has two
para-symmetries~$r_1$, $r_2$, see Remark~\ref{rem.para};
these symmetries generate a certain finite
dihedral group~$\DG{2n}$. Let $c:=r_1r_2$ be the generator of the cyclic
subgroup $\CG{n}\subset\DG{2n}$; it is the rotation through~$2\alpha$,
where $\alpha$ is the angle between the two axes.

\begin{figure}[h]
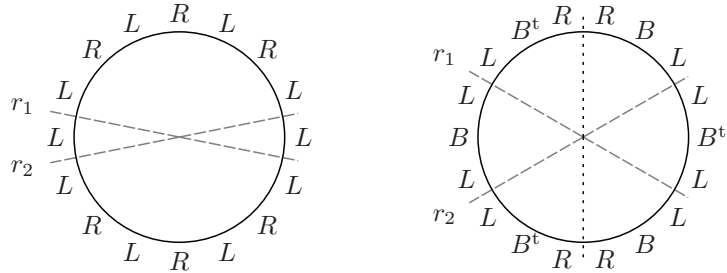

\cpic{Vn}\qquad\qquad\cpic{n-3}
\caption{Diagrams with two para-symmetries}\label{fig.Vn}
\end{figure}

\subsubsection{The two para-symmetries have a common anchor
\rm(see Figure~\ref{fig.Vn}, left)}\label{ss.common.anchor}
In this case,
the $\DG{2n}$-action
on~$\diagram$ is obviously transitive and, starting from
an appropriate anchor,
we arrive at $g\sim V_m$, where $n=2m+1$ and
\[
V_m:=L^2(LR)^mL^2(LR)^m\sim L^2(RL)^mL^2(RL)^m,\quad m\ge0.
\label{eq.Vm}
\]
In particular, $n$ is odd. It is immediate that $\diagram$ has no other
para-symmetries, as it has only four pairs of consecutive occurrences
of~$L$, which could serve as anchors.

\subsubsection{The two para-symmetries have no common anchors
\rm(see Figure~\ref{fig.Vn}, right)}\label{ss.no.anchors}
Consider the orbits of the $\CG{n}$-action
on~$\diagram$. Call an orbit \emph{special} or \emph{ordinary} if it,
respectively, does or does not contain an anchor.
Each ordinary
orbit is `constant', \ie, is either $L^n$ or $R^n$.
To analyze a special orbit, start with an anchor~$a$ of~$r_1$ and observe
that $c$ preserves the letter $a\ra c^i$ unless $i=0\bmod n$
(in this case, $r_1$ preserves~$a$ and $r_2$ reverses $a\ra r_1$, so that
$a\ra c$ is an~$R$)
or $i=k:=[n/2]\bmod n$. In the latter case, if $n=2k$ is even, then
$a\ra c^k$ is an anchor for~$r_1$; otherwise, if $n=2k+1$ is odd, then
$a\ra c^kr_1$ is an anchor for~$r_2$.

Thus, we conclude that $n=2k+1$ must be odd, as otherwise $a\ra c^k$, which
is an~$R$, would be an anchor for~$r_1$. Furthermore, there are four
special orbits of~$\CG{n}$, each one being of the form $LR^kL^k$
(\emph{in the orbit cyclic order}, which may differ from the cyclic order
restricted from~$\diagram$),
where the
first and the $(k+1)$-st letters are anchors for~$r_1$ and~$r_2$,
respectively.

Assume that there is a third para-symmetry~$r$. Together with~$r_1$
and~$r_2$, it generates a dihedral group $\DG{2m}\supset\DG{2n}$ and, since
$\CG{n}\subset\DG{2m}$ is a normal subgroup, $r$ takes $c$-orbits to
$c$-orbits, reversing their orbit order. Unless $n=3$, a special orbit is
taken to a special one, with one of the two anchors contained in the orbit
preserved and the other elements reversed. If $n=3$, a special orbit $LRL$
can be taken to $L^3$, with the two copies of~$L$ preserved. In both cases,
$r$ shares an anchor with~$r_1$ or~$r_2$ and
$g\sim V_m$ for some~$m$, which is a contradiction.
\qed

\corollary[of the proof]\label{cor.symmetric}
In the case of subsection~\ref{ss.no.anchors}, the union of all special
orbits is symmetric with respect to the two reflections~$s_1$, $s_2$ whose
axes bisect the angles between~$r_1$ and~$r_2$.
\endcorollary

\proof
Indeed, since $s_1cs_1=c':=r_2r_1$, the orbit $\{a_1\ra c^i,i\in\Z\}$
starting from an anchor~$a_1$ of~$r_1$ is taken (with the letters preserved)
to the orbit $\{a_2\ra(r_2r_1)^i,i\in\Z\}$ starting from the anchor
$a_2:=a_1\ra s$ of~$r_2$, and the latter orbit is also special.
\endproof

The next corollary refines the statement of Theorem~\ref{th.<=2}.
\corollary[of the proof] The $2$-factorizations corresponding to two distinct
para-symmetries of the cyclic diagram of a hyperbolic element $g\in\MG$
are not weakly equivalent.
\endcorollary
\begin{proof}
According to Corollary~\ref{cor.A1A2}, the $2$-factorizations are weakly equivalent if and only if the two para-symmetries are isomorphic, \ie, related by a rotation symmetry of the cyclic diagram. Since the axes cannot be orthogonal, see subsection~\ref{ss.no.anchors}, this rotation would give rise to more axes, which would contradict  to Theorem~\ref{th.<=2}.
\end{proof}

\corollary\label{cor.g=h1orh2}
If  a hyperbolic $2$-factorizable element $g$ is a power $h^n$
for some $h \in \Gamma$, then $n=1$ or $2$
and, in the latter case, one has
$g\sim(L^2A)^2$ for a word $A$
in $\{L,R\}$
such that $A\trans=A$.
\endcorollary

\proof
Under the assumptions, in addition to a para-symmetry~$r$,
the diagram $\diagram_g$ has a rotation symmetry $c$ of
order~$n>1$, and hence also para-symmetries $c^{-i}rc^i$, $i=0,\ldots,n-1$.
In view of subsection~\ref{ss.no.anchors}, it follows that $n\le2$, as
otherwise $\diagram_g$ would have two para-symmetries with orthogonal axes
(if $n=4$) or more than two para-symmetries (if $n=3$ or $n\ge5$).
\endproof

\remark
From Corollary~\ref{cor.g=h1orh2}, it follows immediately that for such $g$, the cutting
period cycle is either the minimal period or at worst twice the minimal
period of the continued fraction expansion.
\endremark

\subsection{Proof of Theorem~\ref{th.weak=strong}}\label{proof.weak=strong}
If $g\sim\X$ or $g\sim R^2$, the $2$-factorization of~$g$ is unique up to
strong equivalence, see the beginning of subsection~\ref{proof.<=2}.

Let $g$ be a hyperbolic element, and assume that $g=f\1gf$ for some
$f\in\MG$ {\em that is not a power of~$g$}. Then
both
$f$ and $g$ are powers
of a hyperbolic element $h\in\MG$,
see Lemma~\ref{lem.centralizer},
and,
due to Corollary~\ref{cor.g=h1orh2}, we have $g=h^2\sim(L^2A)^2$ and
$A\trans=A$ (hence $\Y A\Y=A\1$).
Modulo~$g$, we can assume that $f=(L^2A)\1$; then
the
$2$-factorization~\eqref{eq.decomposition} and its conjugate by~$f$ differ by
one Hurwitz move. (Geometrically, one can argue that the skeleton~$\skeleton$
of the
monodromy group has a central symmetry and the two $2$-factorizations are
obtained from two symmetric markings of~$\skeleton$.)

\begin{figure}[h]
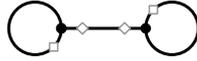

\cpic{L4}
\caption{The skeleton corresponding to $g=L^4$}\label{fig.L4}
\end{figure}

Finally, if $g\sim L^4$, the corresponding skeleton~$\skeleton$ is as shown
in Figure~\ref{fig.L4}. It has four markings with respect to which the
monodromy at infinity is~$L^4$, see the figure, and the corresponding marked
skeletons split into two pairs of isomorphic ones, resulting in two strong
equivalence classes of $2$-factorizations:
\[*
L^4=R\cdot(R\1L^2)R(R\1L^2)\1
=LRL\1\cdot(LR\1L^2)R(LR\1L^2)\1.
\] Note that the two classes are conjugate
by~$L$.
\qed

\subsection{Elements admitting two $2$-factorizations}\label{s.=2}
Let $n=2k+1$ and consider the word
$w_{1/n}:=l(lr)^kl(lr)^k$ in the alphabet~$\{l,r\}$.
Denote by $w[i]$, $i\ge0$,
the $i$-th
letter of a word~$w$, the indexing starting from~$0$.
Pick an odd integer $1\le m<n$ prime to~$n$ and let
$w_q$, $q:=m/n$, be the word in~$\{l,r\}$ of length~$2n$ defined by
\[*
w_q[i]=w_{1/n}[mi\bmod2n],\quad i=0,\ldots,2n-1.
\]
Given a word~$B$ in $\{L,R\}$, let
$w_q\{B\}$ be the word obtained from $w_q$ by inserting a
copy of~$B$ between $w_q[2i]$ and $w_q[2i+1]$
and a copy of~$B\trans$ between $w_q[2i+1]$ and
$w_q[2i+2]$, $i=0,\ldots,n-1$. Finally, let $W_q(B)$ be
the word obtained from $w_q\{B\}$ by the substitution
$l\mapsto L^2$, $r\mapsto R^2$.

\theorem\label{th.=2}
An element $g\in\MG$ admits two distinct strong equivalence classes of
$2$-factorizations if and only if
either $g\sim V_m$,
$m\ge0$, see~\eqref{eq.Vm},
or $g\sim W_q(B)$, where $0<q<1$ is a
rational number with odd numerator and denominator and $B$ is any
word in~$\{L,R\}$, possibly empty.
\endtheorem

\begin{proof}
It has been explained in subsection~\ref{proof.<=2} that each element
$g\sim\X$ or $g\sim R^2$ admits a unique $2$-factorization,
whereas an element
$g\sim V_m$, $m\ge0$, admits two $2$-factorizations
(which are weakly equivalent if $g\sim L^4=V_0$,
see subsection~\ref{proof.weak=strong} for more details on this case).
Thus, it remains to consider a hyperbolic element~$g$
that is not conjugate to any~$V_m$, $m\ge0$.

Consider the cyclic diagram $\diagram=\diagram_g$.
According to subsection~\ref{ss.no.anchors}, the two para-symmetries~$r_1$
and~$r_2$
of~$\diagram$ have no common anchors and their axes are at an
angle~$\alpha$
of the form $\pi m/n$, where $n=2k+1\ge3$ is odd and $m$ is prime to~$n$.
Choosing for~$\alpha$ the minimal positive angle and replacing it, if
necessary, with $\pi-\alpha$, we can assume that $m$ is also odd and
$0<m<n$, so that $q:=m/n$ is as in the statement.
Consider the orbits of the rotation $c:=r_1r_2$.
The union of the special orbits, see subsection~\ref{ss.no.anchors},
is uniquely
determined by the angle~$\alpha$: if $m=1$, then
$g\sim W_{1/n}(\varnothing)$
(with the ordinary orbits disregarded),
\cf. Figure~\ref{fig.Vn}, right; otherwise,
each orbit is `stretched' $m$ times and `wrapped' back around the circle,
so that $g\sim W_q(\varnothing)$.
In the union of the special orbits, adjacent to each semiaxis of each
symmetry contained in~$\DG{2n}$ is a pair of equal letters, either both~$L$
or both~$R$; these pairs are encoded by, respectively, $l$ and~$r$ in the
word $w_{1/n}$ used in the definition of~$W_q$. These pairs divide the circle
into $2n$ arcs, which are occupied by the ordinary orbits and, taking into
account the full $\DG{2n}$-action, one can see that the union of
all ordinary
orbits has the form $B,B\trans,\ldots,B,B\trans$, where $B$ is the portion of this
union in one of the arcs, see Figure~\ref{fig.Vn}, right; it can be any word
in~$\{L,R\}$.
\end{proof}

\subsection{Relation to real structures}\label{ss.relationreal}
Here, we discuss the elements of~$\MG$ that admit both a $2$-factorization
and a real structure.

Clearly, an elliptic element $g\sim\X$ and a parabolic element $g\sim R^2$
have this property. In both cases, the only $2$-factorization is
real. Furthermore, in both cases we have both types of
real structures (or real Lefschetz fibrations, see
Remark~\ref{rem.bm.real.geometric}):
for  $\X=R\cdot L\1$,
the action of~$\hat\tau_1$ preserves the $2$-factorization, whereas that of
$\hat\tau$ with $\tau=\tau_2R\1$
changes it by the Hurwitz move~$\Gs_1$;
for $R^2=R\cdot R$, the action of  $\hat{\tau}_{2}$ can be
regarded as either preserving
the $2$-factorization or changing it by~$\Gs_1$.

A parabolic element $g\sim L^4$ has two strong equivalence classes of
$2$-fac\-tor\-iza\-tions and four real structures, as can be easily seen from its
cyclic diagram. Both $2$-factorizations are real with respect to two of the
real structures and are interchanged by the two others.

\theorem\label{th.2.and.real}
Assume that a hyperbolic element $g\in\MG$  is real and
admits a $2$-factorization~$\BM$. If $\BM$ is real,
then it is unique, and
$g$ has a unique real structure.
Otherwise,
$g$ has two $2$-factorizations, both non real, which are interchanged by the
real structure.
\endtheorem

\begin{proof}
Under the assumptions,
the cyclic diagram $\diagram:=\diagram_g$ has a para-symmetry
(the $2$-factorization)~$r$ and a symmetry (the real structure)~$s$.
Then $r':=srs$ is also a para-symmetry and, unless $r'=r$, the two
$2$-factorizations corresponding to~$r$ and~$r'$
are interchanged by the real structure.

If $r'=r$, \ie, $r$ is real, the axes of~$r$ and~$s$ are orthogonal.
(Since $g\not\sim L^4$, the two axes cannot coincide.) If
there were another para-symmetry $r_1\ne r$, then $r$, $r_1$, and
$r_1':=sr_1s$ would define three distinct $2$-factorizations, which would
contradict to Theorem~\ref{th.<=2}. Similarly, another symmetry $s_1\ne s$
would generate, together with~$s$, a dihedral group~$\DG{2n}$, $n\ge3$,
giving rise to $n$ distinct para-symmetries.
\end{proof}

\remark
The proof of Theorem~\ref{th.2.and.real} gives us a complete characterization
of real hyperbolic elements~$g$ admitting a $2$-factorization~$\BM$.

The $2$-factorization~$\BM$ of~$g$ is real if and only if $g\sim L^2AL^2A\trans$ for a
\emph{palindromic} word~$A$ in $\{L,R\}$.

Otherwise, there are two $2$-factorizations and we have either
$g\sim V_{m}$, $m\ge1$, or $g\sim W_q(B)$, see Theorem~\ref{th.=2}.
In the former case, $g$ has two real structures,
the corresponding symmetries of the
cyclic diagram having orthogonal axes.
In the latter case, due to Corollary~\ref{cor.symmetric},
the union of the special orbits, \ie,
the part $W_q(\varnothing)$, is symmetric with respect to two
reflections~$s_1$, $s_2$ whose axes are distinguished as those
bisecting the `odd' and `even' angle between the axes of the para-symmetries
(respectively, the horizontal and vertical axes in Figure~\ref{fig.Vn}).
The symmetry~$s_1$ is a real structure on $W_q(B)$ if and only if $B$ is
palindromic,
whereas $s_2$ is a real structure if and only if $B=\varnothing$.
(It is worth mentioning that the two non-equivalent factorizations of $V_n$ or
$W_q(B)$ with $B$ palindromic differ by a global conjugation in
the group
$\PGL(2,\Z)$.)
\endremark

\remark
If $\hat\tau(\BM)$
is strongly equivalent to~$\BM$,
then $\hat\tau$ preserves the monodromy group~$G$
of~$\BM$ and; hence,
induces an orientation reversing symmetry of the skeleton
$G\backslash\MG$. Clearly, any such symmetry of a skeleton~$\skeleton$ as in
Proposition~\ref{44} \emph{with at least one Farey branch} must interchange
the two monogons of~$\skeleton$.
Hence, any real $2$-factorization~$\BM$ of a
\emph{hyperbolic} element of~$\MG$ represents
a real Lefschetz fibration with a pair of complex conjugate singular fibers,
see Remark~\ref{rem.bm.real.geometric}; in other words, it is real
in the sense $\hat\tau(\BM)=\BM$ for some real structure~$\tau$.
\endremark

\subsection{Further observations}
For practical purposes the following observation is useful, as it eliminates
most matrices as not admitting a $2$-factorization.

\proposition\label{prop.trace}
If an element $g\in\tMG$
factors into a product of two Dehn twists, then
$(2-\trace g)$ is a perfect square.
\endproposition

\remark
This is definitely not a sufficient condition;
for a counterexample one can take
the element $R^3LR^2=(L^2RL^3)\trans$ of trace~$7$.
\endremark

\proof
Up to conjugation, we can assume that the two Dehn twists
constituting the product are $R=t_{\ba}$ and $A:=t_{[p,q]}$ for some $[p,q]\in\CH$,
$\gcd(p,q)=1$. Since
\[*
A=\bmatrix1-pq&-q^2\\p^2&1+pq\endbmatrix,
\]
one has $\trace RA=2-q^2$; on the other hand,
trace is a class function.
\endproof

As a consequence of the proof, we conclude that,
for each integer~$q$, there does exist an element $g\in\tMG$ of trace $2-q^2$ which is a product of two
Dehn twists in $\tMG$. For an element $g\in\MG$, one should check whether $2\pm\trace g$ is a perfect
square. Proposition~\ref{prop.trace} has a geometric meaning:
the number $2-\trace g$ is the square of the
symplectic product of the eigenvectors of the two Dehn twists.

For
another necessary condition, consider a finite group $G$ and fix an
ordered sequence of conjugacy classes represented by elements
$g_1,\ldots,g_r\in G$. Then the number $N(g_1,\ldots,g_r)$ of solutions to the
equation $x_1\cdot\ldots\cdot x_r=\id$, $x_i\sim g_i$, $i=1,\ldots,r$, is
given by the following Frobenius type formula, see~\cite{Arad.Herzog}:
\[*
N(g_1,\ldots,g_r)=\frac{\ls|g_1|\ldots\ls|g_r|}{\ls|G|}
 \sum\frac{\chi(g_1)\ldots\chi(g_r)}{\chi(\id)^{r-2}},
\]
where $\ls|\,\cdot\,|$ stands for the size of the conjugacy class and the
summation runs over all irreducible characters of~$G$.
Applying this formula to the images of $g_1=g_2=R$, $g_3=g\1$ in a finite
quotient of~$\tMG$, we have the following statement.

\proposition\label{prop.characters}
If an element $g\in\tMG$
factors into a product of two Dehn twists, then, for each positive
integer~$n$, one has
\[*
\sum\chi(R)^2\chi(g\1)\chi(\id)\1\ne0,
\]
the summation running over all irreducible characters of the group
$\SL(2,\Z_n)$.
\done
\endproposition

Note that all irreducible characters of the groups
$\SL(2,\Z_p)$ for $p$ prime are known, see,
\latin{e.g.},~\cite{Naimark.Stern}, and, for each prime~$p$, the condition in
Proposition~\ref{prop.characters} can be
checked effectively
in terms of certain Gauss sums.
At present, we do not know whether an analogue of the Hasse principle holds
for the
$2$-factorization problem, \ie, whether Propositions~\ref{prop.trace}
and~\ref{prop.characters} together constitute a sufficient condition for the
existence of a $2$-factorization.

\section{Real elliptic Lefschetz fibrations}\label{S.Lf}

\subsection{Lefschetz fibrations}\label{s.Lf}

Let $X$ be a compact connected oriented smooth $4$-manifold and
$B$ a compact connected smooth oriented surface. A
\emph{Lefschetz fibration} is a
surjective smooth map
$p\:X\to B$ with the following properties:
\roster*
\item
$p(\partial X)=\partial B$ and the restriction
$p\:\partial X\to\partial B$ is a submersion;
\item
$p$ has  finitely many critical points, which are all in the
interior of~$X$, and all critical values are pairwise distinct;
\item
about each critical point~$x$ of~$p$, there are local charts
$(U,x)\cong(\C^2,0)$ and $(V,b)\cong(\C^1,0)$, $b=p(x)$, in which
$p$ is given by $(z_1,z_2)\mapsto z_1^2+z_2^2$.
\endroster
The restriction of a Lefschetz fibration to the set
$B^\sharp$
of regular values of~$p$ is a locally trivial fibration with all
fibers closed connected oriented surfaces;
the \emph{genus} of~$p$ is the genus of a generic fiber.
Lefschetz fibrations of genus one are called \emph{elliptic}.

An \emph{isomorphism} between
Lefschetz fibrations is a pair of
orientation preserving diffeomorphisms of the total spaces and the
bases commuting with the projections.
The \emph{monodromy} of a Lefschetz
fibration is the monodromy of its restriction to~$B^\sharp$.
As it follows from the local normal form in the definition,
the local monodromy (in the positive direction) about a singular
fiber is the positive Dehn twist about a certain simple closed
curve, well defined up to isotopy; this curve is called the
\emph{vanishing cycle}. The singular fiber
itself is obtained from a close nonsingular one by contracting the
vanishing cycle to a point to form a single node.
A singular fiber is irreducible (remains connected after resolving
the node) if and only if its vanishing cycle is not
null-homologous. If the vanishing cycle bounds a disk, the
singular fiber contains a sphere, which necessarily has
self-intersection $(-1)$, \ie, is a topological analogue of a
$(-1)$-curve. As in the analytic case, such a sphere can be blown
down.
The fibration is called \emph{relatively minimal}
if its singular fibers do not contain $(-1)$-spheres, \ie, none of
the vanishing cycles is null-homotopic.

From now on, we only consider {\em relatively minimal elliptic
Lefschetz fibrations over the sphere $B=S^2$}. After choosing a
base point $b\in B^\sharp$ and fixing an isomorphism
$H_1(p\1(b))=\CH$,
the monodromy of such a fibration becomes a homomorphism
$\pi_1(B^\sharp,b)\to\tMG$, and it is more or less clear
(see~\cite{Moishezon:LNM} for a complete proof) that, up to
isomorphism, the fibration is determined by its monodromy.
By the Riemann--Hurwitz formula, $\chi(X)=r$, where $r$ is the
number of singular fibers.

\theorem[Moishezon, Livn\'{e}~\cite{Moishezon:LNM}]\label{th.Moishezon.Livne}
Up to isomorphism, a relatively minimal elliptic
Lefschetz fibration $X\to S^2$ is determined by
the Euler characteristic $\chi(X)$, which
is subject to the restrictions $\chi(X)\ge0$ and
$\chi(X)=0\bmod12$.
\pni
\endtheorem

Since for any $k\ge0$ there exists an elliptic
surface $E(k)$ with $\chi(E(k))=12k$, it follows that any elliptic
Lefschetz fibration $p\:X\to S^2$ is \emph{algebraic}, \ie, $X$
and~$S^2$ admit analytic structures with respect to which $p$ is a
regular map.

\definition\label{def.Jacobian}
A
\emph{Jacobian Lefschetz fibration} is a relatively minimal
elliptic Lefschetz fibration $p\:X\to B\cong S^2$
equipped with a distinguished section $s\:B\to X$ of~$p$.
Isomorphisms of such fibrations are required to commute with the
sections.
\enddefinition

According to Theorem~\ref{th.Moishezon.Livne}, any elliptic
Lefschetz fibration over~$S^2$ admits a section, which is unique
up to automorphism.

\subsection{Real Lefschetz fibrations}\label{s.Lf.real}

Mimicking algebraic geometry (\cf. subsection~\ref{s.tc.real} below),
define a \emph{real structure} on a
smooth oriented $2d$-manifold~$X$ as an involutive autodiffeomorphism
$c_{X}\:X\to X$ with the following properties:
\roster*
\item
$c_{X}$ is orientation preserving (reversing) if $d$ is even
(respectively, odd);
\item
the \emph{real part} $X_\R:=\Fix c_{X}$ is either empty or of pure
dimension~$d$.
\endroster
A \emph{real Lefschetz fibration} is a Lefschetz fibration
$p\:X\to B$ equipped
with a pair of real structures $c_X\:X\to X$ and
$c_B\:B\to B$ commuting with~$p$. Such a fibration is
\emph{totally real} if all its singular fibers are real.
(Auto-)homeomorphisms of real Lefschetz fibrations are supposed to
commute with the real structures. A Jacobian Lefschetz fibration
is real if the distinguished section is real, \ie, commutes with
the real structures.

Recall that for any real structure~$c$ on~$X$ one has the
\emph{Thom--Smith inequality}
\[
\Gb_*(X_\R)\le\Gb_*(X),
\label{eq.Thom-Smith}
\]
where $\Gb_*$ stands for the total Betti number with
$\Z_2$-coefficients.
If \eqref{eq.Thom-Smith} turns into an equality, the real
structure (or the real manifold~$X$) is called \emph{maximal}.
If $X$ is a closed surface of genus~$g$, we have
$\Gb_0(X_\R)\le g+1$.

From now on, we assume that the base~$B$ is the sphere~$S^2$ and
the real part~$B_\R$ is a circle~$S^1$, \ie, $c_B$ is maximal;
sometimes, $B_\R$ is referred to as
the \emph{equator}. A real Lefschetz fibration equipped with a
distinguished orientation of~$B_\R$ is said to be \emph{directed};
a \emph{directed \rom(auto-\rom)homeomorphism} of such fibrations
is an \hbox{(auto-)}\penalty0 homeomorphism preserving the distinguished
orientations. The fibers over $B_{\R}$ inherit real structures from $c_{X}$;
they are called
\emph{real fibers}.

A large supply of real Jacobian Lefschetz fibrations is provided
by real Jacobian elliptic surfaces, see subsection~\ref{s.tc-es}.
Such fibrations are called \emph{algebraic}; formally, a
real (Jacobian) Lefschetz fibration $p\:X\to B$ is
\emph{algebraic} if $X$
and~$B$ admit analytic structures with respect to which $p$
(and~$s$) are
holomorphic and $c_X$, $c_B$ are anti-holomorphic. It turns out
that some (in a sense, most) Lefschetz fibration are \emph{not}
algebraic; the realizability of a given fibration by an elliptic
surface is one of the principal questions addressed in this
paper, see subsection~\ref{s.algebraic}.

\subsection{Necklace diagrams}\label{s.necklace}

Define a \emph{broken necklace diagram}
as a
nonempty
word in the \emph{stone} alphabet $\{\Circ,\Square,>,<\}$.
Associate to each stone its \emph{dual} and \emph{inverse} stones
and its \emph{monodromy}
(an element of~$\MG$)
as shown in Table~\ref{tab.stones}.
\Table[h]
\caption{Necklace stones}\label{tab.stones}
\def\rel{\joinrel\Rightarrow\joinrel\joinrel\Relbar\joinrel}
\hbox to\hsize{\hss\vbox{\halign{%
 &\tabstrut\hss\quad$#$\quad\hss\cr\toprule
\text{Segment}&\text{Stone}&\text{Dual}&\text{Inverse}&\text{Monodromy}\cr
\midrule
{\circ}{\rel}{\circ}&\Circ&\Square&\Circ&\Y\X^2\Y\X^2\Y\cr
{\times}{\joinrel\rel\joinrel}{\times}&\Square&\Circ&\Square&\X^2\Y\X^2\cr
{\!\times}{\joinrel\rel}{\circ}&{>}&{<}&{<}&\X\Y\cr
{\circ}{\rel\joinrel}{\times}&{<}&{>}&{>}&\Y\X\cr
\bottomrule\crcr}}\hss}
\endTable
Then, given a broken necklace diagram~$\necklace$, we can define
its
\roster*
\item
\emph{monodromy} $\bm(\necklace)\in\MG$, which is
obtained by replacing each stone with
its monodromy and evaluating the resulting word in~$\MG$,
\item
\emph{dual diagram}~$\necklace^*$, obtained by replacing each stone with
its dual, and
\item
\emph{inverse diagram}~$\necklace\1$, obtained by replacing each stone
with its inverse and reversing the order of the stones.
\endroster
Note that the operations of dual and inverse commute
with each other
and that for any
diagram~$\necklace$ one has
$\bm(\necklace^*)=\Y\cdot\bm(\necklace)\cdot\Y$ and $\bm(\necklace\1)=\hat\tau_{1}(\bm(\necklace))$.
Furthermore, the symmetric group~$\SG{n}$ acts on the set
$\BND(n)$ of broken necklace diagrams of length~$n$.
For any \emph{cyclic} permutation $\Gs\in\SG{n}$ one has
$(\necklace\ra\Gs)^*=\necklace^*\ra\Gs$ and $(\necklace\ra\Gs)\1=\necklace\1\ra\Gs\1$; thus, on
$\BND(n)$ there is a well defined action of the group
$\CG2\times\DG{2n}$ generated by the dual, inverse, and cyclic
permutations.

\definition\label{def.ND}
An \emph{oriented necklace diagram} $\necklace$ is an element of the
quotient set
$\BND(n)/\CG{n}$ by the subgroup $\CG{n}$ of cyclic
permutations or, equivalently, a cyclic word in the stone
alphabet.
A \emph{\rom(non-oriented\rom) necklace diagram} is
an element of the quotient $\BND(n)/\DG{2n}$ by the subgroup
generated by the cyclic permutations and the inverse.
\enddefinition

With real trigonal curves in mind, define also oriented
\emph{flat} and \emph{twisted}
necklace diagrams as elements of the quotients
$\BND(n)/\CG2\times\CG{n}$ and
$\BND(n)/\CG2\times\tilde{\Bbb Z}_{2n}$, respectively. Here,
$\CG2$ acts \latin{via} $\necklace\mapsto \necklace^*$,
$\CG{n}$ is the subgroup of
cyclic permutation, and $\tilde{\Bbb Z}_{2n}$ acts \latin{via}
the \emph{twisted shifts}
$S_1S_2\ldots S_n\mapsto S_2\ldots S_nS_1^*$.
In both cases, the non-oriented versions are defined by further
identifying the orbits of~$\necklace$ and $\necklace\1$.

\begin{figure}[h]
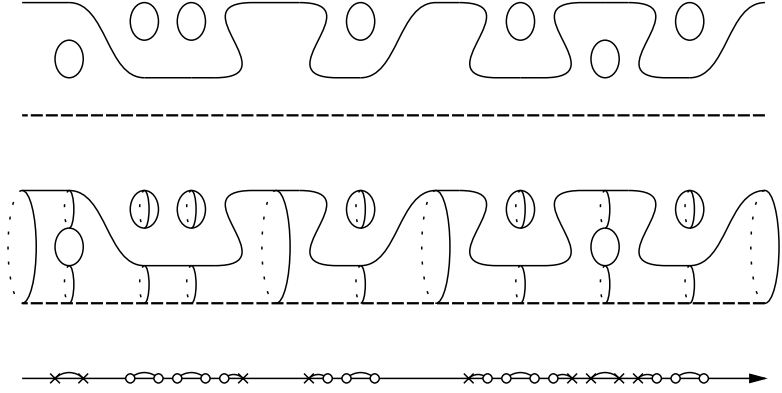

\centerline{\cpic{03-topology}}
\caption{A non-hyperbolic trigonal curve (top),
a covering Jacobian
surface (middle), and its uncoated necklace
diagram (bottom); the horizontal dotted lines represent
the distinguished sections}\label{fig.topology}
\end{figure}

Consider a directed Jacobian Lefschetz fibration $p\:X\to B$ and
assume that it has at least one real singular fiber. The
restriction $p_\R\:X_\R\to B_\R$ can be regarded as an
$S^1$-valued Morse function, and one can assign an index~$0$, $1$,
or~$2$ to each real singular fiber, \latin{i.q\.} critical point
of~$p_\R$. The real part of each real nonsingular fiber is
nonempty (as there is a section); hence it consists of one or two
circles, see~\eqref{eq.Thom-Smith}, and the number of circles
alternates at each singular fiber.
Define the \emph{uncoated necklace diagram} of~$p$ as the
following decoration of the oriented circle~$B_\R$:
\roster*
\item
each singular fiber of index~$0$ or~$2$ is marked with a
\white-,
and each singular fiber of index~$1$ is marked with a \cross-;
\item
each segment connecting two consecutive singular fibers over which
nonsingular fibers have two real components is doubled.
\endroster
A typical real part~$X_\R$ and its
uncoated necklace diagram are shown in
Figure~\ref{fig.topology}, middle and bottom, respectively.

\definition\label{def.Lf->ND}
The \emph{oriented necklace diagram} $\necklace(p)$ of
a directed Jacobian Lefschetz fibration $p\:X\to B$ is the cyclic
word in the stone alphabet obtained by replacing each double
segment of its uncoated necklace diagram with a single stone as
shown in Table~\ref{tab.stones}. In the presence of a
base point~$b$ inside one of the simple segments of~$B_\R$, one can
also speak about the \emph{broken necklace diagram} $\necklace_b(p)$
of~$p$, with
the convention that the first stone~$S_1$
is the immediate successor of~$b$.
\enddefinition

For example, the necklace diagram of the
fibration
shown in Figure~\ref{fig.topology} is
\[*
\ND{*..<>.>.<*>.}\,.
\]
(In~\cite{Salepci:necklace}, necklace diagrams are drawn in the
oriented circle~$B_\R$, and we respect this convention by drawing
a `broken' necklace. For long diagrams we will also use the
obvious multiplicative notation for associative words.)
According to the following theorem, a
\emph{totally real} fibration is uniquely recovered from its
necklace diagram.

\theorem[see~\cite{Salepci:necklace,Salepci:thesis}]\label{th.Nermin}
Given $k>0$,
the map $p\mapsto \necklace(p)$
establishes a bijection between the set of
isomorphism classes of \rom(directed\rom)
totally real Jacobian Lefschetz fibrations with $12k$ singular
fibers
and the set of \rom(oriented\rom) necklace diagrams of length~$6k$
and monodromy~$\id\in\MG$.
\pni
\endtheorem

The classification of totally real Lefschetz fibrations for the small values
of~$k$ is also found in~\cite{Salepci:necklace,Salepci:thesis}. For $k=1$,
there are $25$ undirected isomorphism classes, among which four are maximal.
For $k=2$, the number of classes is~$8421$.

\subsection{Generalizations}\label{s.Lf.general}
Let~$\necklace$ be a broken necklace diagram. A
\emph{$w$-pendant}
on~$\necklace$ is a \emph{strong} Hurwitz equivalence class
of $w$-factorizations~$\BM$ of~$\bm(\necklace)$.
The $(\CG2\times\DG{2n})$-action on the set $\BND(n)$ is
extended to pairs $(\necklace,\BM)$ as follows:
\roster*
\item
the inverse $(\necklace,\BM)\1$ is $(\necklace\1,\hat\tau_{1}(\BM))$;
\item
the dual $(\necklace,\BM)^*$ is $(\necklace^*,\Y \BM\Y)$;

\item
the cyclic permutation $1\mapsto2\mapsto\ldots$ acts \latin{via}
$\necklace=S_1\ldots S_n\mapsto S_2\ldots S_nS_1$
and $\BM\mapsto P_1\1\BM P_1$, where $P_1$ is the monodromy
of~$S_1$.
\endroster
An
\emph{oriented $w$-pendant necklace diagram}
is an orbit of the
cyclic permutation action on the set of pairs $(\necklace,\BM)$ as above;
a (non-oriented) $w$-pendant necklace diagram is obtained by the
further identification of the orbits of $(\necklace,\BM)$ and
$(\necklace,\BM)\1$. The \emph{length} of a $w$-pendant necklace diagram
represented by $(\necklace,\BM)$ is the
length~$\mathopen|\necklace\mathclose|$, the number of stones on $\necklace$.

\remark\label{rem.twisted}
An \emph{oriented flat $w$-pendant necklace diagram} is defined as an orbit of
the further action $(\necklace,\BM)\mapsto(\necklace,\BM)^*$.
In the case of twisted necklace diagrams, both the monodromy
and the notion of $w$-pendant should be defined slightly differently. Namely,
given a broken necklace diagram
$\necklace$, let $\tbm(\necklace):=\bm(\necklace)\Y$.
The twisted shift by the cyclic permutation
$\Gs\:1\mapsto2\mapsto\ldots$ acts \latin{via}
$\tbm(\necklace\ra\Gs)=P_1\1\tbm(\necklace)P_1$, and we can define a
\emph{twisted $w$-pendant} as a strong equivalence class of
$w$-factorizations~$\BM$
of $\tbm(\necklace)$. The twisted action of $\CG2\times\tilde{\Bbb Z}_{2n}$
extends to pairs $(\necklace,\BM)$ in the same way as above, and an
\emph{oriented twisted $w$-pendant necklace diagram}
is defined as an orbit set of this action.
The non-oriented analogues are defined as above, by the further
identification of the orbits of $(\necklace,\BM)$ and $(\necklace,\BM)\1$.
\endremark

\begin{figure}[h]
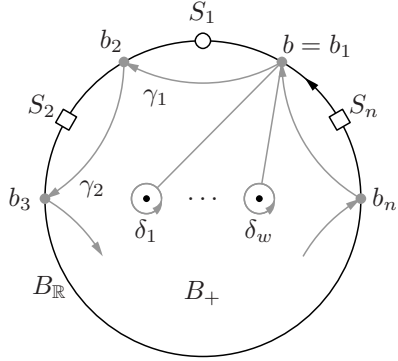

\centerline{\cpic{03-monodromy}}
\caption{The monodromy of a real Lefschetz fibration}\label{fig.monodromy}
\end{figure}

Let $p\:X\to B$ be a directed Jacobian Lefschetz fibration with
$r>0$ real and $w\ge0$ pairs of complex conjugate singular fibers.
Denote by $B_+\subset B$ the closed hemisphere inducing the chosen
orientation of the equator~$B_\R$. Decorate~$B_\R$ as explained in
subsection~\ref{s.necklace} and remove from~$B_+$ the union of
some
disjoint regular neighborhoods of the stones, \latin{i.q\.} double
segments; denote the resulting closed disk by~$\disk$ and let
$\disk^\sharp=\disk\cap B^\sharp$. Choose a base point
$b\in\disk\cap B_\R$ and pick a geometric basis
$\{\Gd_1,\ldots,\Gd_w\}$ for the group $\pi_1(\disk^\sharp,b)$,
see Figure~\ref{fig.monodromy} (where black dots denote non-real
singular fibers).

The real structure~$c :=c_{X}|_{F_{b}}$
in the real fiber~$F_b$ over~$b$
is conjugate to~$\tau_1$; it
gives rise to
a distinguished pair of opposite bases $\pm(\ba,\bb)$ in the
homology $H_1(F_b)$,
which
are defined by the condition that
$\ba\pm\bb$ should be a $(\pm1)$-eigenvector of~$c_*$. Thus, there
is a canonical, up to sign, identification $H_1(F_b)=\CH$ and the
monodromies $\bm(\Gd_i)$ project to well defined elements
$\bm_i\in\MG$, $i=1,\ldots,w$. Let
$\BM_b(p)=(\bm_1,\ldots,\bm_w)$.

\lemma\label{lem.D-BM}
The
strong
equivalence class of the $w$-factorization
$\BM_b(p)$ is indeed a
$w$-pendant on the broken necklace diagram $\necklace_b(p)$. A change of
the base point~$b$ used in the definition results in a cyclic
permutation action on the pair $(\necklace_b,\BM_b)$.
\endlemma

\proof
According to~\cite{Salepci:thesis}, the monodromy~$P_i$ of a
stone~$S_i$ is the $\MG$-valued monodromy along a path~$\Gg_i$
connecting two points~$b_i$ and~$b_{i+1}$, right before and right
after~$S_i$, and circumventing~$S_i$ in the clockwise direction,
see Figure~\ref{fig.monodromy}.
(To obtain a well defined element
of~$\MG$, in the fibers over both points one should use the
canonical bases described above.) Hence, the first statement of
the lemma follows from the obvious relation
$\Gg_1\cdot\ldots\cdot\Gg_n\sim[\partial\disk]$. For the second
statement, it suffices to notice that, changing the base point from
$b=b_1$ to~$b_2$, one can take for a new geometric basis for
$\pi_1(\disk^\sharp,b_2)$
the set $\{\Gg_1\1\Gd_i\Gg_1\}$, $i=1,\ldots,w$.
\endproof

\theorem\label{th.Lf}
The map sending $p\:X\to B$ to the class of the pair
$(\necklace_b(p),\BM_b(p))$ establishes a bijection between the set of
isomorphism classes of \rom(directed\rom) Jacobian Lefschetz
fibrations with
$2n>0$ real and $w$ pairs of complex conjugate singular fibers and the set of \rom(oriented\rom)
$w$-pendant necklace diagrams of length $n$.
\endtheorem

\proof
Due to Lemma~\ref{lem.D-BM}, the map in question is well defined,
and to complete the proof it suffices to show that
a Lefschetz fibration can be recovered from a pair $(\necklace,\BM)$
uniquely up to isomorphism.
The necklace diagram~$\necklace$ gives rise to
a unique, up to isomorphism, totally real directed Jacobian
Lefschetz fibration over an equivariant regular neighborhood~$U$
of the equator~$B_\R$ (see~\cite{Salepci:thesis} for details; this
statement is an essential part of the proof of
Theorem~\ref{th.Nermin}). The complement $B\sminus U$ consists of
two connected components $B^\circ_\pm$, and $\BM$ is a
$w$-factorization of the monodromy
$\bm(\partial B^\circ_+)=\bm(\necklace)$;
due to~\cite{Moishezon:LNM}, this factorization determines a
unique extension of the fibration from $\partial B^\circ_+$
to~$B^\circ_+$. The extension to the other half~$B^\circ_-$ is
defined by symmetry.
\endproof

\remark
It
is not easy to decide whether a given necklace diagram $\necklace$ admits
a $w$-pendant.
There are simple criteria
for $w=0$
(one must have $\bm(\necklace)=\id$),
$w=1$ ($\bm(\necklace)$ must be a Dehn twist), and $w=2$ (the criterion is
given by Theorem~\ref{th.>=1}).
In general, one can lift $\bm(\necklace)$ to a degree~$w$ element
in the braid group~$\BG3$
and apply S.~Orevkov's
quasipositivity criterion~\cite{Orevkov:quasipositivity}:
a $w$-pendant exists if and only if the lift is quasipositive.
A lift
of degree~$w$ exists (and then is unique) if and only if
$\deg\bm(\necklace)=w\bmod6$,
where $\deg\:\MG\onto\CG6$ is the abelianization epimorphism, with the
convention that $\deg R=1$. Obviously, this condition is necessary for the
existence of a $w$-pendant.
\endremark

\remark
If $w=0$ or~$1$,
a necklace diagram~$\necklace$ obviously admits at most one
$w$-pendant. If $w=2$, there are at most two $w$-pendants, see
Theorems~\ref{th.<=2} and~\ref{th.=2}. It follows that at most two
isomorphism classes of
real Jacobian Lefschetz fibrations with two pairs of complex conjugate
singular fibers may share the same necklace diagram (equivalently, fibered
topology of the real part).
\endremark

We used {\tt Maple} to compute the numbers of undirected isomorphism classes
of real Jacobian Lefschetz fibrations for some small values of~$k$ and~$w$
(where $12k$ is the total number of singular fibers and $w$ is the number of
pairs of complex conjugate ones). For $k=1$, the numbers are $25$ ($w=0$),
$28$ ($w=1$), and $24$ ($w=2$); for $k=2$, they are $8421$ ($w=0$)
and $15602$ ($w=1$).
(For $k=w=2$, the computation is too long.)
In all examples,
a fibration with $w>0$ pairs of conjugate singular fibers can be obtained
from one with $(w-1)$ pairs by converting the pair of real fibers
constituting an arrow type stone to a pair of conjugate ones. We do not know
how general this phenomenon is.

\subsection{Counts}\label{s.counts}

We conclude this section with a few simple counts. Let $p\:X\to B$
be a real Jacobian Lefschetz fibration, $\chi(X)=12k>0$, and let
$\necklace=\necklace(p)$. We assume that $\necklace\ne\varnothing$, so that
$0<\ls|\necklace|\le6k$. Denote by $\ncirc:=\ncirc(\necklace)$
the number of stones of type~$\diamond$, $\diamond \in\{\Circ,\Square,>,<\}$.
Then one has
\[
\Gb_*(X)=\chi(X)=12k
\label{eq.X}
\]
and
\[
\Gb_*(X_\R)=2(\nCirc+\nSquare)+4,\qquad
\chi(X_\R)=2(\nCirc-\nSquare),
\label{eq.XR}
\]
see~\cite{Salepci:necklace,Salepci:thesis}
or Figure~\ref{fig.topology}. In particular,
$\nCirc+\nSquare\le6k-2$, see~\eqref{eq.Thom-Smith}, and $X$ is
maximal if and only if $\nCirc+\nSquare=6k-2$.
Alternatively, $X$ is maximal if and only if
\[
\NO<+\NO>+w=2,
\label{eq.maximal}
\]
where $w\le2$ is the number of pairs of complex
conjugate singular fibers.

\section{Real trigonal curves}\label{S.tc}

\subsection{Trigonal curves}\label{s.tc}

A \emph{Hirzebruch surface}~$\Sigma_d$, $d>0$, is a
geometrically ruled rational surface
with a (unique) exceptional section~$E$ of self-intersection~$-d$.
We denote the ruling by $p\:\Sigma_d\to B\cong\Cp1$; its fibers
are called the fibers of~$\Sigma_d$.
A \emph{\rom(proper\rom)
trigonal curve} is a reduced curve $C\subset\Sigma_d$
disjoint from the exceptional section and intersecting
each fiber at three points (counted with multiplicities); in other
words, the restriction $p\:C\to B$ is a map of degree three.
A \emph{singular fiber} of
a trigonal curve is a fiber of~$\Sigma_d$ intersecting the curve
geometrically at fewer than three points; equivalently, singular
are the
fibers over the critical values of the restriction $p|_C$.
In this paper, we consider \emph{almost generic} trigonal curves
only, \ie, we assume that all critical points of $p|_C$ are
simple. Such a curve is nonsingular and irreducible;
hence it has genus $g(C)=3d-2$ (the adjunction formula)
and the number of singular fibers is $6d$ (the Riemann--Hurwitz
formula).

Given a point~$b\in B$, let $F_b$ be the fiber $p\1(b)$ and let
$F^\circ_b$ be the \emph{affine fiber} $F_b\sminus E$.
It is an affine complex line. Hence, in the presence of a trigonal
curve, one can speak about the \emph{zero section} $Z$ sending
each point $b\in B$ to the barycenter of the three points
$C\cap F^\circ_b$, and about the fiberwise convex hull
$\conv C\subset\Sigma_d\sminus E$.

\subsection{Real trigonal curves}\label{s.tc.real}

Recall that a \emph{real structure} on an algebraic (analytic)
variety~$X$ is an anti-holomorphic involution $c_{X}\:X\to X$. A pair
$(X,c_{X})$ is called a \emph{real algebraic variety}; usually, the
real structure~$c$ is understood and we speak about a real
algebraic variety~$X$.
Given a real structure, the fixed point set $X_\R:=\Fix c_{X}$ is
called the \emph{real part} of~$X$.
A maximal real algebraic variety, see~\eqref{eq.Thom-Smith}, is
usually called an \emph{$M$-variety}.

Up to isomorphism, a Hirzebruch surface~$\Sigma_d$ admits a unique
real structure $c_d$ with nonempty real part $(\Sigma_d)_\R$,
the latter being a torus or a Klein bottle for $d$ even or odd,
respectively. This real structure~$c_d$ descents to a certain real
structure~$c_{B}$ on the base~$B$, so that the ruling~$p$ is real.
In what follows, when speaking about a real Hirzebruch surface,
we assume such a pair $(c_d,c_{B})$ fixed. Any real automorphism
$\Gf\:\Sigma_d\to\Sigma_d$ induces a $c_{B}$-equivariant
autohomeomorphism (in fact, a real automorphism) $\Gf_B\:B\to B$.
Such an autohomeomorphism~$\Gf_B$ (and the original
automorphism~$\Gf$) is said to be \emph{directed} if it preserves
the orientation of~$B_\R$.
More generally, a \emph{directed
autohomeomorphism} of a real Hirzebruch surface
$p\:\Sigma_d\to B$ is
an orientation preserving $c_d$-equivariant fiberwise
autohomeomorphism whose descent to~$B$ preserves
the orientations of~$B$ and~$B_\R$, \cf.
subsection~\ref{s.Lf.real}.

A trigonal
curve $C\subset\Sigma_d$ is \emph{real} if it is $c_d$-invariant;
then, the restriction of~$c_d$ to~$C$ is a real structure on~$C$.
In affine coordinates, such a curve is given by a polynomial with
real coefficients. By a \emph{deformation} of real trigonal curves
we mean an equivariant deformation (a path in the space of real
polynomials) in the class of almost generic curves. Two curves are
said to be
\emph{\rom(directedly\rom) deformation equivalent}
if they differ by a
deformation and/or (directed) real automorphism
of~$\Sigma_d$.

Fix a real trigonal curve $C\subset\Sigma_d$ and consider the
restriction $p_\R\:C_\R\to B_\R$. If each fiber of~$p_\R$ consists
of three points, the curve~$C$ is called \emph{hyperbolic}. In
this case, the real part~$C_\R$ consists of three (if $d$ is even)
or two (if $d$ is odd) components; in the former case, each
component is mapped onto~$B_\R$ homeomorphically; in the latter
case, one `central' component is mapped homeomorphically and the
other is a double covering. It can be shown that all hyperbolic
curves (in a given surface) are deformation equivalent, see,
\latin{e.g.},~\cite{DIK:elliptic}.
If $C$ is not hyperbolic, its real part looks like shown in
Figure~\ref{fig.topology}, top.
More precisely, $C_\R$ has one \emph{long
component} that is mapped onto~$B_\R$ and, possibly, a number of
\emph{ovals}, necessarily unnested. The long component may contain
some \emph{zigzags} (the Z-shaped fragments in the figure).
For a formal definition, consider a maximal, with respect to
inclusion, segment $I\subset B_\R$ with the
property that each point $b\in I$ has at least two pull-backs
under~$p_\R$. If the pull-back $p_\R\1(I)$ is disconnected, one of
its components is an oval; otherwise, the pull-back is called a
zigzag.
With a certain abuse of the language, the projections of the ovals
and zigzags to~$B_\R$ (\ie, maximal segments $I\subset B_\R$
as above)
are also referred to as ovals and zigzags,
respectively, \cf. subsection~\ref{s.topology} below.

\subsection{The covering elliptic surface}\label{s.tc-es}

Let $C\subset\Sigma_d$ be a real trigonal curve. If $d$ is even,
the double covering $X\to\Sigma_d$ ramified at $C+E$ is a
Jacobian
elliptic surface ($E$ being the section), and the real
structure~$c_d$ lifts to two \emph{opposite} real structures
$c_\pm\:X\to X$ that differ by the deck translation of the
covering.
Disregarding the analytic structure, one can consider
the corresponding real varieties~$X_\pm$ as real Jacobian
Lefschetz fibrations.

The necklace diagrams $\necklace(X_\pm)$ are dual
to each other, and each of them determines the embedded topology
of the real part $C_\R\subset(\Sigma_d)_\R$. Hence, the latter can
be encoded by the pair $\necklace(X_\pm)$, \ie, by a flat necklace
diagram.

If $d$ is odd, a covering elliptic surface only exists over the
complement of the fiber $F_b$ over a point $b\in B$. In this case,
choosing for~$b$ a generic real point with one preimage
under~$p_\R$ and analyzing the dependence of~$X_\R$ on~$b$, one
can see that the real part $C_\R\subset(\Sigma_d)_\R$ is
encoded by a twisted necklace diagram.

In both cases, the ovals of~$C_\R$ correspond to the $\Square$ and
$\Circ$ type stones, whereas the zigzags correspond to the arrow
type stones, see Figure~\ref{fig.topology}.

\subsection{Dessins}\label{s.dessins}

In appropriate (real) affine coordinates $(x,y)$, a
(real) trigonal curve
$C\subset\Sigma_d$ can be given by its \emph{Weierstra{\ss}
equation}
\[
y^3+g_2(x)y+g_3(x)=0,
\label{eq.Weierstrass}
\]
where $g_2$ and~$g_3$ are some (real)
polynomials in~$x$ of degree at
most~$2d$ and~$3d$, respectively, and the discriminant
$\Delta(x):=-4g_2^3-17g_3^2$ is not identically zero. The
\emph{$j$-invariant} of~$C$ is defined as the
meromorphic function
\[
j_C\:B\to\Cp1=C\cup\{\infty\},\quad j_C:=-\frac{4g_2^3}\Delta.
\label{eq.j}
\]
If $C$ is real, so is $j_C$.

If $C$ is almost generic, by a small
equisingular deformation it can be made
\emph{generic}, \ie, such that the $j$-invariant~$j_C$ has
\emph{generic branching behavior} in the sense
of~\cite{Moishezon:LNM}; the latter means that all zeroes of~$j_C$
are triple, all zeroes of $(j_C-1)$ are double, all poles of $j_C$
are simple, and all critical values
of~$j_C$ other than~$0$ or~$1$ are also simple. (We emphasize that
these properties are highly \emph{non}-generic for a map
$B\to\Cp1$, but they do correspond to truly generic trigonal
curves.)
For generic real trigonal
curves we define
\emph{\rom(directed\rom) strict deformation
equivalence} as the equivalence relation generated by
the (directed) real automorphisms of~$\Sigma_d$ and equivariant
equisingular deformations \emph{in the class of generic curves}.

Fix a generic trigonal curve~$C$ and define its \emph{dessin}
$\dessin:=\Dssn C$ as the embedded graph $j_C\1(\Rp1)\subset B$
decorated as follows:
\roster*
\item
the pull-backs of~$0$, $1$, and~$\infty$ are \black--,
\white--, and \cross-vertices, respectively;
\item
an edge is colored solid, bold, or dotted if its image belongs to
$(-\infty,0)$, $(0,1)$, or $(1,\infty)$, respectively and
is directed
according to the canonical orientation of~$\Rp1$, \ie, the
standard linear order on~$\R$;
\item
the critical points of~$j_C$ with
real critical values distinct from~$0$, $1$, or~$\infty$
are considered \emph{monochrome} vertices, respectively
solid, bold, or dotted.
\endroster
If $C$ is real, its dessin~$\dessin$ is invariant with respect to
the real structure~$c$ on~$B$ and the real part $B_\R$ is a union
of edges and vertices of~$\dessin$. The properties of the graph
thus obtained are summarized in the following definition.

\definition\label{def.dessin}
Let~$B$ be the sphere~$S^2$ and $c\:B\to B$ the reflection against
the equator.
A \emph{\rom(real\rom) dessin} is a $c$-invariant embedded directed
graph $\dessin\subset B$
decorated with the
following additional structures (referred to as the
\emph{colorings}
of the edges and vertices of~$\dessin$, respectively):
\roster*
\item
each edge of~$\dessin$ is of one of the three kinds: solid, bold,
or dotted;
\item
each vertex of~$\dessin$ is of one of the four kinds:
\black-, \white-, \cross-, or monochrome
(the vertices of the first three kinds being called
\emph{essential}),
\endroster
and satisfying the following conditions:
\roster
\item\label{tg-boundary}
the equator of~$B$ is a union of edges and vertices of~$\dessin$,
and each monochrome vertex is at the equator;
\item\label{tg-valency}
each \black--, \white--, \cross--, or monochrome vertex has
valency~$6$, $4$, $2$, or~$4$, respectively;
\item\label{tg-oriented}
the orientations of the edges of~$\dessin$ induce an orientation
of the boundary of the complement $B\sminus\dessin$;
\item\label{tg-monochrome}
all edges incident to a monochrome vertex are of the same kind;
\item\label{tg-cross}
\cross-vertices are incident to incoming dotted edges and
outgoing solid edges;
\item\label{tg-black}
\black-vertices are incident to incoming solid edges and
outgoing bold edges;
\item\label{tg-white}
\white-vertices are incident to incoming bold edges and
outgoing dotted edges;
\item\label{tg-cycles}
$\dessin$ has no directed \emph{monochrome cycles}, \ie,
directed cycles with all
edges of the same kind and all vertices monochrome.
\endroster
In items~\ref{tg-cross}--\ref{tg-white}, the lists are complete,
\ie, vertices cannot be incident to edges of other kinds or with a
different orientation.

The equator $\Fix c$ is called the \emph{real part} of the dessin;
the edges and vertices of~$\dessin$ that are in the equator are
called \emph{real}, whereas the other vertices are called
\emph{inner}.

Two dessins are said to be
\emph{\rom(directedly\rom) homeomorphic} if
they are related by a (directed)
orientation preserving $c$-equivariant
autohomeomorphism of~$B$. Note that directedly homeomorphic
dessins are, in fact, equivariantly isotopic.
\enddefinition

If $C$ is a generic curve in~$\Sigma_d$, then $\deg j_C=6d$ and,
hence, the numbers of \black--, \white--, and \cross-vertices of
the dessin $\Dssn C$ are $2d$, $3d$, and $6d$, respectively. The
number $3d$ of \white-vertices is called the \emph{degree} of
the curve~$C$ and dessin $\Dssn C$.

\begin{conv.}
In view of the symmetry,
in the figures we only draw the portion
of a dessin contained in one of the two hemispheres, which is
represented by a disk. The real part of the dessin
(the boundary of the disk) is
shown by a thick grey line, \cf., for example,
Figure~\ref{fig.cubics}
on page~\pageref{fig.cubics}.
When speaking about directed
dessins,
we choose the closed hemisphere $B_+$ whose
orientation induces the fixed orientation of the equator $B_\R$.
\end{conv.}

\theorem[see~\cite{DIK:elliptic,Orevkov:Riemann}]\label{th.dessin}
The map $C\mapsto\Dssn C$ establishes a bijection between the set
of \rom(directed\rom) strict deformation equivalence classes of
generic real trigonal curves and that of
\rom(directed\rom) homeomorphism classes of dessins.
\pni
\endtheorem

\remark
The notion of strict deformation equivalence is not very
meaningful from the topological point of view, as some codimension
one
degenerations of the $j$-invariant do not affect the topology of
the curve. It is shown in~\cite{DIK:elliptic} that deformation
equivalence classes of almost generic curves are in a one-to-one
correspondence with certain equivalence classes of dessins, where
two dessins are considered equivalent if they are related by a
sequence of homeomorphisms and certain elementary moves. We omit
the description of these moves as they are not essential in the
case of $M$-curves, which are only considered in this paper.
\endremark

\definition\label{def.skeleton}
The union of the bold edges and \black-- and \white-vertices of
the dessin $\Dssn C$ is called the \emph{skeleton} of~$C$ and is
denoted by $\Sk C$.
If $\Dssn C$ has no bold monochrome vertices,
$\Sk C$ is a
regular skeleton in the sense of
the definition given in subsection~\ref{ss.subgroups}.
\enddefinition

\subsection{Topology in terms of dessins}\label{s.topology}

The topology of the real part $C_\R\subset(\Sigma_d)_\R$ is easily
recovered from the dessin $\dessin:=\Dssn C$. It is immediate
from~\eqref{eq.Weierstrass}, \eqref{eq.j}, and the definition
of~$\dessin$ that
\roster
\item\label{dssn-cross}
the singular fibers of~$C$ are those over the \cross-vertices
of~$\dessin$,
\item\label{dssn-white}
the points of intersection $C\cap Z$ are
over the \white-vertices of~$\dessin$, and
\item\label{dssn-dotted}
a point~$b$ inside a real edge~$e$ of $\dessin$ has three
preimages under $p_\R$ if and only if $e$ is dotted.
\endroster
It follows that the ovals and zigzags of~$C_\R$, regarded as
subsets of~$B_\R$, are the maximal dotted segments in~$B_\R$; any
such segment is bounded by two \cross-vertices and is allowed to
contain a number of \white-- or monochrome vertices inside.
In view of item~\ref{dssn-white} above, a maximal dotted segment
in~$B_\R$ is an oval (zigzag) if and only if it contains an even
(respectively, odd) number of \white-vertices.

The uncoated necklace diagram~$N$ of the covering elliptic
surface~$X$
of $C$ is also recovered from~$\dessin$. (If the degree of~$C$ is
odd, we should fix a base point~$b$ inside a solid or bold real
edge of~$\dessin$ and speak about the uncoated diagram broken
at~$b$.)
The \cross-- and \white-points of~$N$ are the \cross-vertices
of~$\dessin$,
and the double segments of~$N$ are the maximal dotted segments
of~$\dessin$.
(At this point, we need to apologize for the notation clash. Unfortunately,
both notation sets are quite well established and it seems unwise to change
them.)
Since $X$ has two opposite real structures, we
cannot distinguish between the \cross-- and \white-type critical
points, but we can compare pairs of points: two critical
points are of the same type (both \cross-\ or both \white-) if and
only if they are separated by an even number of \white-vertices
of~$\dessin$. All assertions are simple consequences of elementary
Morse theory; they are obvious from
Figure~\ref{fig.topology}.

Given a necklace diagram~$\necklace$, denote by $\NO{\fam0 ess}$ the
number of \emph{simple} segments of the corresponding uncoated
diagram connecting pairs of critical points of opposite types,
\ie, those of the form
${\circ}{\joinrel\rightarrow\joinrel\joinrel\relbar\joinrel}\!\!{\times}$
or
${\times}\!\!{\joinrel\rightarrow\joinrel\joinrel\relbar\joinrel}{\circ}$.
(Such segments are called \emph{essential}.)
Since any real Jacobian elliptic surface~$X$ is the double
covering of the Hirzebruch surface $\Sigma_d=X/{\pm\id}$ ramified
at~$E$ and a certain real trigonal curve~$C$,
we have the following necessary condition for a necklace diagram~$\necklace$
to be algebraic.

\proposition[see~\cite{Salepci:necklace,Salepci:thesis}]\label{th.essential}
Let $X\to B$ be a real Jacobian elliptic surface
with $\chi(X)=12k>0$.
Then its
necklace diagram~$\necklace$ is subject to the inequalities
\[*
\NO{\fam0 ess}\le2k,\qquad
\NO{\fam0 ess}+\NO<+\NO>\le6k.
\]
\endproposition

\proof
The
second statement follows from the fact that each zigzag (an
arrow type stone) and each essential segment contains an odd
number, hence at least one, of \white-vertices of the dessin, and
$\nwhite=6k$. For the first one, observe that a \white-vertex
inside an essential segment is separated from
each of the two \cross-vertices
bounding the segment by at least one \black-vertex (as the type of
the edges must change from bold to solid, \cf.
Figure~\ref{fig.cubics} on page~\pageref{fig.cubics}),
and $\nblack=4k$.
\endproof

\subsection{The monodromy}\label{s.tc-monodromy}

Let $\disk\subset B$ be a closed disk. A continuous section
$s\:\disk\to\Sigma_d$ of~$p$ is called \emph{proper} (with respect
to a fixed trigonal curve~$C$) if its image is disjoint from
both~$E$ and $\conv C$. Since
the disk~$\disk$ is contractible and all
fibers $F_b\sminus(E\cup\conv C)$ are connected, a proper
section exists and
is unique up to homotopy in the class of proper sections.

Fix a trigonal curve~$C$, a disk~$\disk$, and a proper
section~$s$.
Assume that the boundary $\partial\disk$ contains no singular
fibers of~$C$ and denote by~$\disk^\sharp$ the disk~$\disk$ with
all singular fibers removed. The restriction
$p\:p\1(\disk^\sharp)\sminus(C\cup E)\to\disk^\sharp$ is a locally
trivial fibration, and one can consider the associated bundle with
the discrete fibers $\Aut\pi_1(F_x\sminus(C\cup E),s(x))$,
$x\in\disk^\sharp$.
This bundle is a covering and, fixing a base point
$b\in\disk^\sharp$ and lifting loops starting from the identity
over~$b$, we obtain a homomorphism
$\tilde\bm\:\pi_1(\disk^\sharp,b)\to\Aut\pi_b$, where
$\pi_b:=\pi_1(F_b\sminus(C\cup E),s(b))$. This homomorphism is called the
\emph{monodromy}; since the
section~$s$ is proper, it actually takes values in the braid group
$\BG3\subset\Aut\pi_b$, where $\pi_b$ is identified with
the free group~$\FG3$ by
means of a geometric basis.

In the sequel, we `downgrade' the monodromy to the modular
group~$\MG$ and consider the composition
$\bm\:\pi_1(\disk^\sharp,b)\to\BG3\onto\MG$,
called the \emph{reduced monodromy}.
If $d$ is even, $\bm$ coincides with the $\MG$-valued reduction of
the monodromy (homological invariant) of the covering elliptic
surface.
The following statement is essentially contained
in~\cite{degt:kplets}, although the conventions and notation
in~\cite{degt:kplets} differ slightly from those accepted in this
paper.

\theorem[see~\cite{degt:kplets}]\label{th.monodromy}
Consider a connected component~$\skeleton_0$ of the intersection
$\Sk C\cap\disk$, and assume that the base point~$b$ is in an edge~$e$
of~$\skeleton_0$. Then, under an appropriate choice of a geometric
basis in the reference fiber~$F_b$,
the diagram
\[*
\CD
\pi_1(\skeleton_0,e)@>i_*>>\pi_1(\disk^\sharp,b)\\
@Vj_*VV@VV\bm V\\
\pi_1(\Sk C,e)@>\val>>\MG
\endCD
\]
commutes, where $i\:\skeleton_0\into\disk^\sharp$
and $j\:\skeleton_0\into\Sk C$
are the inclusions.
\pni
\endtheorem

Now, assume that the curve~$C$ is real. Orient the real
part~$B_\R$ and consider the positive hemisphere~$B_+$.
Choose some disjoint regular neighborhoods $U_i\subset B$ of the
singular fibers of~$C$ and let
$\disk=B_+\sminus\bigcup U_i$, the union running
\emph{over the real singular fibers only} (\cf. similar
construction in subsection~\ref{s.Lf.general} and
Figure~\ref{fig.monodromy}).
Pick a base point~$b$ in the boundary $\partial\disk$,
make the other necessary choices,
and consider the reduced monodromy
$\bm\:\pi_1(\disk^\sharp,b)\to\MG$.

\definition\label{def.BMG}
The image $\BMG(C):=\Im\bm\subset\MG$ is called the
\emph{monodromy group} of~$C$; the element
$\bminf:=\bm[\partial\disk]\in\BMG(C)$
is called the \emph{monodromy at infinity}.
\enddefinition

The following statement is straightforward.

\proposition\label{prop.BMG}
The pair
$(\BMG,\bminf)$,
$\bminf\in\BMG\subset\MG$,
is determined by the
curve~$C$
and orientation of~$B_\R$ up to conjugation.
The conjugacy class of
the pair
$(\BMG,\bminf)$
is invariant under
directed autohomeomorphisms of the pair
$(\Sigma_d,C)$\rom; in particular, it is a directed
deformation invariant of real trigonal curves.
\done
\endproposition

The conjugacy class of the monodromy at infinity depends on the
real part~$C_\R$ only. Indeed, choose the base point~$b$ real and outside
the zigzags and ovals and represent the real part by a broken
(at~$b$) necklace diagram~$\necklace$,
flat or twisted,
see subsection~\ref{s.tc-es}. Then,
up to conjugation, one has
$\bminf=\bm(\necklace)$ if $d$ is even and
$\bminf=\tbm(\necklace)$
if $d$ is odd.
(Note that the conjugacy class of $\tbm(\necklace)$ is preserved by the
twisted shifts used in the definition of twisted
diagrams, see Remark~\ref{rem.twisted}.)

\subsection{Real trigonal $M$-curves}\label{s.M-curves}

In view of~\eqref{eq.Thom-Smith}, the real part~$C_\R$ of
an $M$-curve $C\subset\Sigma_d$ has $3d-1$ connected components.
Hence, unless $d=1$, such a curve is non-hyperbolic, its real
part has $3d-2$ ovals, and, since each oval and each zigzag
consume two real singular fibers, one has $z+w=2$, where $z$ is
the number of zigzags and $w$ is the number of pairs of complex
conjugate singular fibers;
in particular, $z,w\le2$.
Comparing
this
to~\eqref{eq.maximal},
we conclude that $C$ is an $M$-curve if and only if any/both
covering elliptic surfaces $X_\pm$ are $M$-varieties. (This
statement trivially holds for hyperbolic curves as well.)

\begin{figure}[h]
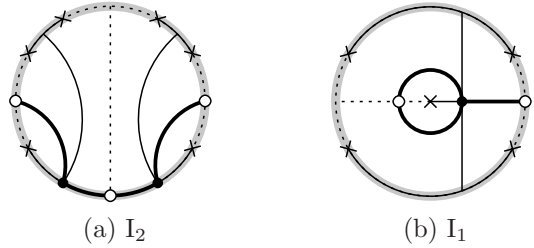

\centerline{\vbox{\halign{&\qquad\hss#\hss\qquad\cr
\cpic{04-cI2}&\cpic{04-cI1}\cr
\noalign{\medskip}
(a) $\I_2$&(b) $\I_1$\cr}}}
\caption{Dessins of $M$-cubics}\label{fig.cubics}
\end{figure}

For $d=1$ (essentially, plane cubics), there are four deformation
families of $M$-curves: one hyperbolic and three non-hyperbolic,
denoted by $\I_z$, $z=0,1,2$:
a curve of type~$\I_z$ has one
oval and $z$ zigzags, see Figure~\ref{fig.cubics} for $z=1$
and~$2$.

To describe the other $M$-curves, we need the operation of
\emph{junction}.
Consider two directed dessins $\dessin^i\subset B^i$, $i=1,2$.
Choose a pair of zigzags $Z^i$ of~$\dessin^i$ with a single
\white-vertex~$v^i$ in each, and let
$I^i\subset\partial B^i_+$ be a segment contained
in the interior of~$Z^i$ and containing~$v^i$ inside.
Pick a homeomorphism $\Gf\:I^1\to I^2$ and consider the
connected boundary sum $B_+:=B_+^1\sqcup_\Gf B_+^2$ and the graph
\[*
\dessin_+:=(\dessin^1\cap B^1_+)\sqcup_\Gf
 (\dessin^2\cap B^2_+)\subset B_+.
\]
Finally, double~$B_+$ to form a new sphere~$B\cong S^2$ and
double~$\dessin_+$ to form a graph $\dessin\subset B$; the two
real \white-vertices~$v^1$, $v^2$ are replaced with a pair of
complex conjugate \white-vertices, and the common endpoints
of~$I^1$, $I^2$ become monochrome vertices of~$\dessin$.
The resulting graph $\dessin\subset B$ is a dessin; it is called
the \emph{junction} of~$\dessin^1$ and~$\dessin^2$ along the pair
of zigzags~$Z^1$, $Z^2$. Up to isotopy, the junction depends only
on the pair of dessins $\dessin^1$, $\dessin^2$, pair of
zigzags~$Z^1$, $Z^2$, and whether
the homeomorphism~$\Gf$ is orientation preserving or
reversing. The zigzags~$Z^1$, $Z^2$ are `consumed' by the
junction, being replaced by a pair of ovals. It follows that an
iterated junction of several dessins does not depend on the order
of the individual operations:
one can start with the disjoint union of all dessins involved and identify
all pairs of segments (which are all disjoint) simultaneously.

An example of junction is shown in
Figure~\ref{fig.junction} on page~\pageref{fig.junction},
where only some essential parts of the dessin are drawn.

\theorem[see~\cite{DIK:elliptic,Orevkov:Riemann}]\label{th.M-curves}
Each directed deformation class of real trigonal $M$-curves
$C\subset\Sigma_d$, $d\ge2$, contains a representative whose
dessin is an iterated junction of $d$ copies of cubic
dessins~$\I_2$, $\I_1$ shown in Figure~\ref{fig.cubics}.
Conversely, any such junction is the dessin of an $M$-curve.
Two \rom(directed\rom) deformation classes are equal if and only
if their corresponding dessins as above
are \rom(directedly\rom)
homeomorphic.
\pni
\endtheorem

According to the number of zigzags,
a dessin of type~$\I_1$ participates in exactly
one junction
(consuming its only zigzag),
whereas a dessin of type~$\I_2$ can participate in
one or two junctions. Furthermore, a dessin of type~$\I_1$ has a
`horizontal' axis of symmetry
preserving the zigzag,
see Figure~\ref{fig.cubics}, whereas a dessin of type~$\I_2$ does
not.
It follows that, after the junction as in
Theorem~\ref{th.M-curves}, the individual blocks form a linear
chain, which
can be encoded by a word in the alphabet
$\{\Mup,\Mdown,\Mside\}$,
with the convention that $\Mside$ can only appear as the first
and/or last letter. Here, $\Mside$ represents a dessin of
type~$\I_1$, and $\Mup$ and $\Mdown$ represent a dessin of
type~$\I_2$, oriented, respectively, as shown in
Figure~\ref{fig.cubics} or upside down.
(For an example, see Figure~\ref{fig.junction} on
page~\pageref{fig.junction}.)
On the set~$M_d$
of such words of length~$d$,
there is an action of the group $\CG2\times\CG2$ generated by the
\emph{vertical flip}~$v$, reversing the order of the letters, and
the \emph{horizontal flip}~$h$, interchanging $\Mup$ and~$\Mdown$.
Each flip is realized by a homeomorphism of dessins
reversing the orientation of the equator.
In these terms, one can restate
Theorem~\ref{th.M-curves} as follows.

\corollary\label{cor.M-curves}
The directed \rom(undirected\rom)
deformation classes of real trigonal
$M$-curves in~$\Sigma_d$ are
in a natural one-to-one correspondence with the orbit set
$M_d/v h$
\rom(respectively, the orbit set
$M_d/\<v,h\>$\rom).
\done
\endcorollary

The number of $\Mside$ type letters in the word representing an
$M$-curve~$C$ equals the number $w\le2$
of pairs of complex conjugate singular fibers of~$C$.

As explained at the beginning of this section,
an $M$-curve~$C$ has at most two zigzags. If
$C$ has at least one zigzag,
the appropriate (flat or twisted)
necklace diagram of~$C$
determines the representation of the dessin $\Dssn C$ in the form
of iterated junction and, hence, the deformation class of~$C$,
see~\cite{DIK:elliptic} for details. If $C$ has no zigzags, this
assertion is no longer true.
The shortest example is the pair of degree~$12$
curves represented by
$\Mcurve{.ud.}$ and $\Mcurve{.du.}$\,: they share the same
oriented flat necklace diagram $\Circ^5\Square^5$, but are not
related by a directed deformation equivalence. The two degree~$30$
curves represented by
\[
\Mcurve{.dudduudu.}\qquad\text{and}\qquad
\Mcurve{.duududdu.}
\label{eq.2.M.curves}
\]
share the same oriented diagram
$(\Circ^5\Square^5\Circ\Square^3)^2$, but are not deformation
equivalent,
directedly or not. An explanation of this phenomenon is given in
Remark~\ref{rem.2.M.curves}
below.
In more details, we address this (non-)uniqueness
question in the next section, see Corollaries~\ref{cor.M.<=2} and~\ref{cor.M.=2}.

\section{Geometric applications}\label{S.appl}

\subsection{Real trigonal $M$-curves are quasi-simple}\label{s.M-simple}
Recall that a deformation family of complex algebraic varieties is called
\emph{quasi-simple}
if, within this family, the equivariant deformation class of a
real variety is determined by the topology of its real structure.  The first
geometric application of our algebraic results is the quasi-simplicity of
real trigonal $M$-curves.

\theorem\label{th.def=dif}
Two
real trigonal $M$-curves $C_1,C_2\subset\Sigma_d$ are
in the same \rom(directed\rom) deformation class if and
only if the
quadruples $(\Sigma_d,C_i,p,c_d)$,
$i=1,2$,
are related
by a \rom(directed\rom) homeomorphism.
\endtheorem

The case where each curve has at least one zigzag is settled
in~\cite{DIK:elliptic}; in this case, non-equivalent curves differ by their
real parts (more precisely, appropriate necklace diagrams,
flat or twisted). Thus, we need to
consider curves without zigzags only (equivalently, those with two pairs of
conjugate singular fibers), and for such curves Theorem~\ref{th.def=dif}
follows from a much stronger statement, Theorem~\ref{th.def=dif.z=0} below.

\definition
A word~$w$ in the alphabet $\{L,R\}$ is called \emph{even}
if all letters occur in~$w$ in pairs or, equivalently, if $w$ can be
represented as a word in $\{L^2,R^2\}$. A parabolic or hyperbolic element
$g\in\MG$ is \emph{even} if, up to conjugation, $g$ is represented by an even
word in $\{L,R\}$. For a hyperbolic element, this condition is equivalent to
the requirement that all entries of the cutting period cycle of~$g$ should be
even.
Finally,
a regular pseudo-tree~$\skeleton$ with two loops, \cf.
Proposition~\ref{44}, is \emph{even} if its monodromy at
infinity is even. Informally, $\skeleton$ is even if its Farey branches
pointing upwards/downwards appear in pairs, see
Remark~\ref{rem.AB}. \enddefinition

\theorem\label{th.def=dif.z=0}
Two real trigonal $M$-curves $C'\subset\Sigma_{d'}$ and
$C''\subset\Sigma_{d''}$ without zigzags
are directedly deformation equivalent \rom(in
particular, $d'=d''$\rom) if and only if the monodromy groups
$\BMG(C')$ and $\BMG(C'')$ are conjugate in~$\MG$.
Furthermore, a subgroup $G\subset\MG$ is the monodromy group of a real
trigonal $M$-curve without zigzags if and only if the skeleton
$G\backslash\MG$
is an even regular pseudo-tree with two loops.
\endtheorem

\proof
Let $C$ be a trigonal curve as in the statement.
According to Theorem~\ref{th.M-curves}, we can assume that the dessin of~$C$
is a junction of two cubic dessins~$\I_1$ and several cubic dessins~$\I_2$,
see Figure~\ref{fig.junction}.
\begin{figure}[h]
\cpic{junction}
\caption{A junction of six $M$-cubics ($\protect\Mcurve{.uudu.}$)}\label{fig.junction}
\end{figure}
(To simplify the figure, we only show
the real part~$B_\R$, maximal real dotted segments, which are all ovals,
junctions,
and the portion of the skeleton $\Sk C$
that is in~$B_+$.) We extend the graph $\Sk C\cap B_+$
to a regular skeleton~$\skeleton'$
by attaching, at
each real \black-vertex, a Farey branch reaching beyond the boundary
to~$B_-$ (shown in dotted bold lines in Figure~\ref{fig.junction}).
Then, the following statements are straightforward:
\roster*
\item
$\skeleton'$ is an even regular pseudo-tree with two loops, \cf.
Proposition~\ref{44};
\item
any even regular pseudo-tree with two loops can be obtained in this way,
starting from a certain junction of $M$-cubics;
\item
the monodromy group $\BMG(C)$ is the stabilizer
of~$\skeleton'$
(see Theorem~\ref{th.monodromy});
\item
the dessin $\Dssn C$ is uniquely recovered from~$\skeleton'$.
\endroster
The first two statements imply the last assertion of the theorem.
The fact that a curve is uniquely determined by its monodromy group follows
from the last two statements and Theorem~\ref{th.M-curves}.
\endproof

In view of our previous results concerning $2$-factorizations, we have the
following corollaries of Theorem~\ref{th.def=dif.z=0}.

\corollary[\cf. Theorem~\ref{th.>=1}]\label{cor.M.>=1}
A necklace diagram~$\necklace$ \rom(flat or twisted\rom)
without arrow type stones is the
diagram of a real trigonal $M$-curve if and only if the monodromy $\bm(\necklace)$
has the form $L^2AL^2A\trans$ for some \emph{even} word~$A$ in $\{L,R\}$.
\done
\endcorollary

\corollary[\cf. Theorem~\ref{th.<=2}]\label{cor.M.<=2}
A necklace diagram~$\necklace$ as in Corollary~\ref{cor.M.>=1} is the diagram of at
most two, up to equivalence, real trigonal $M$-curves. In other words, at
most two equivalence classes of curves may have homeomorphic real parts.
\done
\endcorollary

\corollary[\cf. Theorem~\ref{th.<=2}]\label{cor.M.=2}
A necklace diagram~$\necklace$
as in Corollary~\ref{cor.M.>=1}
gives rise to two equivalence classes of real trigonal
$M$-curves if and only if the monodromy $\bm(\necklace)$
has the form $W_q(B)$ for some \emph{even} word~$B$ in $\{L,R\}$.
\done
\endcorollary

\remark\label{rem.2.M.curves}
The
two $M$-curves given by~\eqref{eq.2.M.curves} which share the same
topology of the real part (\latin{i.q.} flat necklace diagram)
correspond to the two
non-equivalent $2$-factorizations of the element $W_{1/3}(L^2R^2)$,
see subsection~\ref{s.=2} and Theorem~\ref{th.=2}.
\endremark

\subsection{Maximal Lefschetz fibrations are algebraic}\label{s.algebraic}
As stated in the introduction, one of the major questions in the theory of
real elliptic Lefschetz fibrations is whether a given fibration can be
realized by an algebraic one. Unlike the complex case, there do exist
non-algebraic real Lefschetz fibrations. Thus, amongst the $25$
undirected
isomorphism classes of totally real fibrations with twelve singular fibers
only seventeen are algebraic,
see~\cite{Salepci:necklace,Salepci:thesis}; the
eight
others are ruled out
by Proposition~\ref{th.essential}. Out of the $8421$ classes of
totally real fibrations
with $24$ singular fibers, at least $4825$ classes are non-algebraic as they
violate Proposition~\ref{th.essential}. At present, we do not know if all
$3596$ remaining classes are algebraic, nor do we know any simple
criterion that would establish that a given fibration \emph{is} algebraic.
(In the case of twelve singular fibers, an analytic structure is constructed
in~\cite{Salepci:necklace,Salepci:thesis} by finding a dessin with the
desired flat necklace diagram.)

In this section, we prove Theorem~\ref{th.algebraic}, closing the question
for maximal fibration.

Consider a
fibration $p\:X\to B$ as in the statement of the theorem and let $\necklace:=\necklace(p)$
be its necklace diagram. Let $w$ be the number of pairs of complex conjugate
singular fibers of~$p$; we have $w\le2$, see~\eqref{eq.maximal}.
The case $w=0$ is covered by Theorem~\ref{th.Nermin.alg}, and it remains to
consider the cases $w=1,2$.

\proof[Proof of Theorem~\ref{th.algebraic}\rom: the case $w=2$]
The diagram~$\necklace$ has no arrow type stones,
see~\eqref{eq.maximal}; hence, $\necklace$ is of the form
$\Square^{i_1}\Circ^{j_1}\ldots\,\Square^{i_s}\Circ^{j_s}$
(assuming
that $\necklace$ has both $\Circ$
and~$\Square$ type stones and breaking it between a $\Circ$ and a $\Square$).
Then, using Table~\ref{tab.stones}, after cancelations we have
\[*
\bm(\necklace)=(RL^{i_1-1}R)(RL^{j_1-1}R)\ldots(RL^{i_s-1}R)(RL^{j_s-1}R).
\]
Thus, in the
cyclic diagram~$\diagram$ of~$\bm(\necklace)$, the copies of~$R$ appear in pairs.
On the other hand,
since $\bm(\necklace)$ admits a $2$-factorization, see Theorem~\ref{th.Lf},
$\diagram$
has a para-symmetry, see Remark~\ref{rem.para}.
It follows that the copies
of~$L$ also appear in pairs (two pairs of anchors and pairs of~$L$ symmetric
to those of~$R$),
\ie, $\bm(\necklace)$ is an even element of~$\MG$.
Hence, the monodromy group of any $2$-factorization of~$\bm(\necklace)$ is
as in Theorem~\ref{th.M-curves} and the corresponding flat pendant necklace
diagram is realized by a real trigonal $M$-curve~$C$.
Due to Theorem~\ref{th.Lf},
one of the two opposite real
Jacobian elliptic surfaces ramified at $C+E$
is isomorphic to~$p$.

If
all stones of~$\necklace$ are of the same type, then, since the monodromy of
each $\Circ$ or $\Square$ type stone is conjugate to~$L$, see
Table~\ref{tab.stones}, we have $\bm(\necklace)\sim L^r$, where $r$ is the
length of~$\necklace$. Hence, $r=4$ and $\necklace=\Circ^4$ or $\Square^4$,
see Theorem~\ref{th.>=1}, and, in view
of Theorem~\ref{th.weak=strong}, $\necklace$ admits two non-equivalent
$2$-pendants.
It is immediate that the four real Lefschetz fibrations obtained in this way
are isomorphic to the four rational real elliptic surfaces ramified over the
two zigzag free trigonal $M$-curves of degree six,
namely $\Mcurve{.uu.}$ and $\Mcurve{.ud.}$.
\endproof

To complete the proof, we need a few observations and a lemma.

Consider a real Jacobian elliptic $M$-surface
$p\:X\to B$ with two pairs of conjugate
singular fibers. In the class of real Jacobian elliptic surfaces,
there are exactly four ways to collide a pair of conjugate singular fibers to
a single real fiber~$F$ of type~$\tilde{\bold{A}}_0^{**}$
(Kodaira's type~$\II$) and perturb~$F$
to produce a pair of real singular fibers (see, \latin{e.g.},
\cite{DIK:elliptic}; in
the realm of real trigonal $M$-curves, this procedure corresponds to
replacing one of the $\Mside$ type letters at an end of the word
representing the junction of $M$-cubics with a
$\Mup$ or $\Mdown$). In two cases, the new necklace diagram has an extra
stone of type $>$, and in the two other cases, it has an extra stone of
type~$<$. Cutting the original diagram $\necklace(p)$ at a base point~$b$ right before
the new stone, we conclude that the pendant necklace diagram of $p$ has two
representatives of the form $(\necklace',(R,g'))$
\rom(for some $g'\in\MG$\rom) and
two representatives of the form $(\necklace'',(L\1,g''))$ \rom(for some
$g''\in\MG$\rom).

Conversely, given a representative of a pendant necklace diagram of the
form, \latin{e.g.}, $(\necklace,(R,\ldots))$, the corresponding real Lefschetz
fibration can be modified,
\emph{in the topological category}, so that a pair of conjugate
singular fibers disappears to produce a $>$ type stone. (Topologically, we
merely remove an equivariant disk surrounding the two fibers and replace it
with another disk, with two real singular fibers and
the same monodromy at the boundary,
see the discussion of the two real structures
on the $2$-factorization $\X=R\cdot L\1$ in subsection~\ref{ss.relationreal}.)
The following lemma asserts that, in the case of a \emph{maximal} Lefschetz
fibration, any such topological modification is one of the two described
above and, hence, can be realized in the algebraic category.

\lemma\label{lem.two.places}
The pendant necklace diagram of a
real Jacobian elliptic $M$-surface
with two pairs of complex conjugate real fibers
has
exactly two representatives of the form $(\necklace',(R,g'))$ \rom(for some
$g'\in\MG$\rom) and
exactly two representatives of the form $(\necklace'',(L\1,g''))$ \rom(for some
$g''\in\MG$\rom).
\endlemma

\begin{figure}[h]
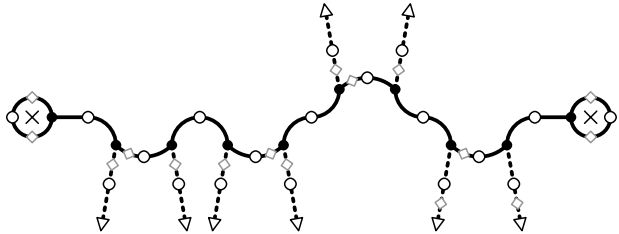

\cpic{points}
\caption{Necklace base points in the skeleton}\label{fig.points}
\end{figure}

\proof
We consider the representatives of the form $(\necklace',(R,g'))$, which
result in the
$>$ type stones. According to the discussion above, two such
representatives do exist; they correspond to the two
algebraic modifications of the fibration. Fix
one of these representatives.
Since the action of~$\Gs_1^2$
is the conjugation by $Rg'=\bm(\necklace')$, the $2$-factorization
of any
other representative in question
is either $(P\1RP,\ldots)$ or $(P\1g'P,\ldots)$, where
$P$ is a monodromy (not necessarily shortest, \ie, possibly multiplied by a
power of $\bm(\necklace')$) from the original base point~$b'$ to the new base point
$b''$. We assert that $P\1RP=R$ if and only if $P=\id$. Indeed, consider the
part $\Sk C\cap B_+$ of the skeleton of the
corresponding real trigonal curve
and extend it to a pseudo-tree~$\skeleton$ as explained in the proof of
Theorem~\ref{th.def=dif.z=0}, \cf. Figure~\ref{fig.junction}. Assume that the
original base point~$b'$ is contained in the boundary of the leftmost cubic
dessin~$\I_1$ in the figure and assign to~$b'$ the base point $e':=e\ra\Y$
of~$\skeleton$, where $e$ is the edge constituting the leftmost monogonal
region.
Starting from~$e'$, `project'
any other base point~$b''$ to~$\skeleton$ by
assigning to~$b''$ the edge $e'':=e'\ra P$. Taking into account
Table~\ref{tab.stones} and using induction, one can easily show that the
edges obtained are either as shown in Figure~\ref{fig.points}
(if the
stone~$S$ preceding~$b'$ is of type~$\Square$) or those in the figure shifted
by~$\Y$ (if $S$ is of type~$\Circ$).
Crucial is the fact that all these points are pairwise distinct. On the other
hand, $e'\ra R^n=e'$ for any $n\in\Z$. Hence,
unless $P=\id$, the monodromy~$P$ is not a power of~$R$ and
we have $P\1RP\ne R$.

It follows, in particular, that the $>$ type stone that may appear in the
cubic dessin~$\I_1$ at the other end of the junction results in a
$2$-factorization of the other form, \ie, $(P\1g'P,\ldots)$,
and the same argument as
above shows that this representative is also unique.
\endproof

\proof[Proof of Theorem~\ref{th.algebraic}\rom: the case $w=1$]
According to~\eqref{eq.maximal}, the
necklace
diagram $\necklace:=\necklace(p)$
has a single arrow type stone~$S$, which we can assume
of type~$>$.
(Otherwise, switch to the dual diagram~$\necklace^*$.)
Replace the two real singular fibers contained in~$S$ with a pair
of conjugate singular fibers (the inverse of the operation described prior
to Lemma~\ref{lem.two.places}).
The new Jacobian Lefschetz fibration is maximal and, according to the first
part of the proof, it is algebraic. Now, due to Lemma~\ref{lem.two.places},
there are only two ways to revert the operation and recreate a $>$ type
stone; both result in algebraic Lefschetz fibrations.
\endproof

As a consequence of Theorem~\ref{th.algebraic}, we have a more effective
description of the deformation classes of maximal real
Jacobian Lefschetz fibrations:
the directed (undirected) deformation classes
of such fibrations are
in a natural one-to-one correspondence with the orbit set
$M_d/v h$
\rom(respectively, the orbit set
$M_d/\<v,h\>$\rom), \cf. Corollary~\ref{cor.M-curves}.

{
\let\.\DOTaccent
\def\cprime{$'$}
\bibliographystyle{amsplain}
\bibliography{degt}

\providecommand{\bysame}{\leavevmode\hbox to3em{\hrulefill}\thinspace}
\providecommand{\MR}{\relax\ifhmode\unskip\space\fi MR }
\providecommand{\MRhref}[2]{%
  \href{http://www.ams.org/mathscinet-getitem?mr=#1}{#2}
}
\providecommand{\href}[2]{#2}
\begin{thebibliography}{10}

\bibitem{Arad.Herzog}
Z.~Arad, J.~Stavi, and M.~Herzog, \emph{Powers and products of conjugacy
  classes in groups}, Products of conjugacy classes in groups, Lecture Notes in
  Math., vol. 1112, Springer, Berlin, 1985, pp.~6--51. \MR{783068}

\bibitem{Artin}
E.~Artin, \emph{Theory of braids}, Ann. of Math. (2) \textbf{48} (1947),
  101--126. \MR{0019087 (8,367a)}

\bibitem{Bardakov}
V.~G. Bardakov, \emph{The structure of a group of conjugating automorphisms},
  Algebra Logika \textbf{42} (2003), no.~5, 515--541, 636, English translation:
  Algebra Logic 42 (2003), no. 5, 287--303. \MR{2025714 (2004k:20072)}

\bibitem{Bogomolov}
Fedor Bogomolov and Yuri Tschinkel, \emph{Monodromy of elliptic surfaces},
  Galois groups and fundamental groups, Math. Sci. Res. Inst. Publ., vol.~41,
  Cambridge Univ. Press, Cambridge, 2003, pp.~167--181. \MR{2012216
  (2004j:14045)}

\bibitem{degt:kplets}
Alex Degtyarev, \emph{Zariski {$k$}-plets via dessins d'enfants}, Comment.
  Math. Helv. \textbf{84} (2009), no.~3, 639--671. \MR{2507257 (2010f:14028)}

\bibitem{degt:monodromy}
\bysame, \emph{{H}urwitz equivalence of braid monodromies and extremal elliptic
  surfaces}, Proc. Lond. Math. Soc. (3) \textbf{103} (2011), 1083--1120.

\bibitem{DIK:elliptic}
Alex Degtyarev, Ilia Itenberg, and Viatcheslav Kharlamov, \emph{On deformation
  types of real elliptic surfaces}, Amer. J. Math. \textbf{130} (2008), no.~6,
  1561--1627. \MR{2464028 (2009k:14071)}

\bibitem{Kulikov.Auroux.Shevchishin}
Vik.~S. Kulikov, D.~Oru, and V.~Shevchishin, \emph{Regular homotopy of
  {H}urwitz curves}, Izv. Ross. Akad. Nauk Ser. Mat. \textbf{68} (2004), no.~3,
  91--114, English translation: Izv. Math. 68 (2004), no. 3, 521--542.
  \MR{2069195 (2005f:14039)}

\bibitem{Kulkarni}
Ravi~S. Kulkarni, \emph{An arithmetic-geometric method in the study of the
  subgroups of the modular group}, Amer. J. Math. \textbf{113} (1991), no.~6,
  1053--1133. \MR{1137534 (92i:11046)}

\bibitem{Moishezon:infty}
B.~Moishezon, \emph{The arithmetic of braids and a statement of {C}hisini},
  Geometric topology ({H}aifa, 1992), Contemp. Math., vol. 164, Amer. Math.
  Soc., Providence, RI, 1994, pp.~151--175. \MR{1282761 (95d:20069)}

\bibitem{Moishezon:LNM}
Boris Moishezon, \emph{Complex surfaces and connected sums of complex
  projective planes}, Lecture Notes in Mathematics, Vol. 603, Springer-Verlag,
  Berlin, 1977, With an appendix by R. Livne. \MR{0491730 (58 \#10931)}

\bibitem{Naimark.Stern}
M.~A. Na{\u\i}mark and A.~I. {\v{S}}tern, \emph{Theory of group
  representations}, Grundlehren der Mathematischen Wissenschaften [Fundamental
  Principles of Mathematical Sciences], vol. 246, Springer-Verlag, New York,
  1982, Translated from the Russian by Elizabeth Hewitt, Translation edited by
  Edwin Hewitt. \MR{793377 (86k:22001)}

\bibitem{Orevkov:Riemann}
Stepan~Yu. Orevkov, \emph{Riemann existence theorem and construction of real
  algebraic curves}, Ann. Fac. Sci. Toulouse Math. (6) \textbf{12} (2003),
  no.~4, 517--531. \MR{2060598 (2005d:14084)}

\bibitem{Orevkov:talk}
\bysame, \emph{On braid monodromy monoid}, Talk at MSRI, 2004.

\bibitem{Orevkov:quasipositivity}
\bysame, \emph{Quasipositivity problem for 3-braids}, Turkish J. Math.
  \textbf{28} (2004), no.~1, 89--93. \MR{2056762 (2005d:20063)}

\bibitem{Ozturk.Salepci:Giroux}
Ferit \"{O}zt\"{u}rk and Nerm\.{\i}n Salepc\.{\i}, \emph{Real open books and
  real contact structures}, to appear, \verb+arXiv:1202.5928+.

\bibitem{Salepci:necklace}
Nerm\.{\i}n Salepc\.{\i}, \emph{Classification of totally real elliptic
  {L}efschetz fibrations via necklace diagrams}, to appear,
  \verb+arXiv:1104.1794+.

\bibitem{Salepci:thesis}
\bysame, \emph{Real {L}efschetz fibrations}, Ph.D. thesis, Universit{\'e}
  {L}ouis {P}asteur et IRMA, 2007.

\bibitem{Salepci:real}
\bysame, \emph{Real elements in the mapping class group of {$T^2$}}, Topology
  Appl. \textbf{157} (2010), no.~16, 2580--2590. \MR{2719402}

\bibitem{Series}
Caroline Series, \emph{On coding geodesics with continued fractions}, Ergodic
  theory ({S}em., {L}es {P}lans-sur-{B}ex, 1980) ({F}rench), Monograph.
  Enseign. Math., vol.~29, Univ. Gen\`eve, Geneva, 1981, pp.~67--76. \MR{609896
  (82h:30052)}

\end{thebibliography}
}

\end{document}